\documentclass[10pt]{amsart}
\usepackage{amscd,amsmath,amssymb,amsthm}
\usepackage{hyperref}
\usepackage[all]{xy}
\usepackage{psfrag}

\newtheorem{thm}{Theorem}[section]
\newtheorem{lemma}{Lemma}[section]
\newtheorem{cor}{Corollary}[section]
\newtheorem{prop}{Proposition}[section]
\newtheorem{defn}{Definition}[section]
\newtheorem{rem}{Remark}[section]

\newcommand{\PP}{\mathcal{P}}  
\newcommand{\w}{\omega}        
\newcommand{\J}{\mathbf{J}}    
\newcommand{\D}{\mathbf{D}}    

\newcommand{\Lie}{\mathbf{L}}   
\newcommand{\R}{\mathbb{R}}    
\newcommand{\Ad}{\mathrm{Ad}}  
\newcommand{\ad}{\mathrm{ad}}  
\newcommand{\im}{\mathrm{im}\, }  
\newcommand{\FL}{\mathbb{F}\mathrm{L}}  
\newcommand{\Sl}{\mathbf{S}}   
\newcommand{\OO}{\mathcal{O}}  
\newcommand{\g}{\mathfrak{g}}
\newcommand{\h}{\mathfrak{h}}
\newcommand{\lir}{\mathfrak{r}}
\newcommand{\q}{\mathfrak{q}}
\newcommand{\be}{\begin{equation}}
\newcommand{\ee}{\end{equation}}
\newcommand{\bea}{\begin{eqnarray}}
\newcommand{\eea}{\end{eqnarray}}
\newcommand{\II}{\mathbb{I}} 
\newcommand{\rII}{\hat{\mathbb{I}}_0} 
\newcommand{\Proj}{\mathbb{P}} 
\newcommand{\Exp}{\mathrm{Exp}\,} 
\newcommand{\restr}[1]{\vrule height3ex width.4pt
depth1.4ex\lower1.4ex\hbox{\scriptsize $\,#1$}}
\newcommand{\rrestr}[1]{\vrule height2ex width.4pt
depth0.9ex\lower0.9ex\hbox{\scriptsize $\,#1$}}

\makeatletter \@addtoreset{figure}{section}
\def\thefigure{\thesection.\@arabic\c@figure}
\def\fps@figure{h, t}
\@addtoreset{table}{bsection}
\def\thetable{\thesection.\@arabic\c@table}
\def\fps@table{h, t}
\@addtoreset{equation}{section}

\makeatother

\allowdisplaybreaks
\def\intprod{\mathbin{\hbox to 6pt{%
                    \vrule height0.4pt width5pt depth0pt
                    \kern-.4pt
                    \vrule height6pt width0.4pt depth0pt\hss}}}

\begin{document}

\title{The symplectic normal space of a cotangent-lifted Action}

\author{M. Perlmutter, M. Rodr\'{\i}guez-Olmos and M.E. Sousa-Dias}

\thanks{M.P.: Institute of Fundamental Sciences, Massey
University. New Zealand. m.perlmutter@massey.ac.nz\\ M. R.-O.
(corresponding author): Ecole Polytechnique F\'ed\'erale de
Lausanne (EPFL), Section de
Math\'ematiques. Lausanne, Switzerland. miguel.rodriguez@epfl.ch\\
E. S.-D.: Departamento de Matem\'atica, Instituto Superior
T\'ecnico. Lisboa, Portugal. edias@math.ist.utl.pt}

\maketitle

\begin{abstract}
For the cotangent bundle $T^*Q$ of a smooth Riemannian manifold
acted upon by the lift of a smooth and proper action by isometries
of a Lie group, we characterize the symplectic normal space at any
point. We show that this space splits as the direct sum of the
cotangent bundle of a linear space and a symplectic linear space
coming from reduction of a coadjoint orbit. This characterization
of the symplectic normal space can be expressed solely in terms of
the group action on the base manifold and the coadjoint
representation. Some relevant particular cases
are explored.
\end{abstract}

\section{Introduction}

  One of the most important tools in the study of the local
geometry of Ha\-mil\-to\-nian $G$-spaces (symplectic manifolds
acted upon by a smooth and proper action of a Lie group $G$ with
an equivariant momentum map) is the so-called Symplectic Slice
Theorem obtained in \cite{GuiSt,Marle},
 and generalized in \cite{LeBa}. This result, which
provides a semiglobal model of a neighborhood of a group orbit by
 an equivariant analogue of Darboux charts, has  proved to
be of fundamental importance in the study of not only the
geometrical properties of Hamiltonian $G$-spaces but also of the
qualitative dynamics  of symmetric Hamiltonian systems defined on
them.

A key ingredient of the Symplectic Slice Theorem  is a linear
space $N$ called the symplectic normal space. The  space $N$ has
played an important role in the study of the Marsden-Weinstein
reduced spaces at singular points of the momentum map (see
\cite{LeBa,OrRa2003,SaLe}), notably for proving that such a
reduced space is Whitney stratified, and can be locally modelled
on a symplectic quotient of $N$ relative to the linear action of a
compact Lie group. On the dynamics side, the symplectic normal
space has been extensively used to produce results on the
stability, persistence and bifurcations for singular relative
elements of symmetric Hamiltonian systems (see for instance
\cite{LeSi,MoRoSt1, Montaldi,OrRa99,RobSD,RobWuLa}). The abstract
space $N$ can be realized as a concrete subspace $V$ of the
tangent space to our symplectic manifold at the point under
consideration, and hence in applications, $V$ provides a concrete
way to study the symplectic normal space $N$.

This work is devoted to the construction of models of $N$ and $V$
when the Hamiltonian $G$-space under study is the cotangent bundle
$T^*Q$ of a $G$-manifold $Q$, equipped with the canonical
symplectic form and acted upon by the cotangent lift of the action
of $G$ on $Q$. The main motivation for this study comes from
geometric mechanics where cotangent bundles are universal phase
spaces,  in the sense that most Hamiltonian systems arising from
classical mechanics are formulated on cotangent bundles of some
configuration spaces, or realized as the symplectic or Poisson
reduced spaces of  appropriate cotangent bundles by the action of
some symmetry group. Despite the importance of cotangent bundles
in mechanics and of the symplectic normal space in the study of
singular points of the momentum map, there has not been an
intensive research in this area. Our approach exploits the
specificity of this bundle structure and of the cotangent-lifted
action. Our aim is to reduce the problem as much as possible to
the study of the geometry of $Q$ and of the $G$-action on it,
rather than on the whole symplectic manifold $T^*Q$, thus
providing a computational foundation for concrete applications.
For this, we use the fact that all the information about the
symplectic geometry of $T^*Q$ and its supported Hamiltonian lifted
action is obtainable from geometric data on $Q$ and the $G$-action
on it.

The paper is organized as follows: In Section~\ref{secWA}, we
review the construction of the space $N$ for a general Hamiltonian
$G$-space. In Section~\ref{secisom}, we first construct a
metric-dependent isomorphism between the tangent space at a point
$p_x$ of $T^*Q$ and $T_xQ \oplus T_x^*Q$ and then we identify this
space with $(\mathfrak{r}\oplus {\bf S})\oplus
(\mathfrak{r}^*\oplus {\bf S}^*)$, where $\Sl$ is a linear slice
at $x$ for the $G$-action on $Q$ and $\mathfrak{r}$ is a linear
slice at the identity for the right $G_x$-action on $G$.  In
Section~\ref{secvectorfields} we introduce a family of local
vector fields on $Q$ and we prove a number of technical results
about them that will be computationally important for the main
results. In Section \ref{secinfinitesimal},
Theorem~\ref{infinitesimalaction} gives the infinitesimal
cotangent-lifted action at $p_x$ in $(\mathfrak{r}\oplus {\bf
S})\oplus (\mathfrak{r}^*\oplus {\bf S}^*)$ coordinates. Our first
main result is in Section \ref{secV}, in
Theorem~\ref{generalsymplectic}, where we provide an explicit
choice for the space $V$ isomorphic to the symplectic normal space
$N$ as a subspace of $T_{p_x}(T^*Q)$ at any point $p_x\in T^*Q$.
Our second main result is Corollary~\ref{generalsymplectic2} which
shows that in general the symplectic normal space at a point
$p_x\in T^*Q$ is isomorphic to the direct sum of the cotangent
bundle of a linear space with the symplectic normal space at
$\mu=\J (p_x)$ for the action of $G_x$ on the coadjoint orbit
through $\mu$. This corollary also provides normal forms for its
symplectic form and momentum map. Finally,
Section~\ref{secspecialV} is devoted to some particular cases of
geometric and dynamical relevance for which the normal form for
the symplectic form is particularly simple.

We view the aforementioned results in sections~\ref{secV} and
\ref{secspecialV} as the main results of the paper. We believe
that the results of subsection~\ref{totallyisotropic} can be
applied to the study of the local properties of a new
stratification for the singular reduced space of $T^*Q$ at zero
momentum, found in \cite{PerRoSD}.
 Also, as the results
contained in subsection~\ref{relativeequilibria} are applicable to
every relative equilibrium of a simple mechanical system, it is
expected that they will contribute to the study of the stability
and bifurcations of relative equilibria for systems of this type.
Indeed in the singular context these results are well-suited to
implement the ideas introduced in the regular case in
\cite{SiLeMar}. There, the authors reduce  the degree of
complexity of the stability analysis by expressing the ``Hessian
criterion" on the symplectic normal space only in terms of the
geometry of $Q$ and the $G$-action on it. This research direction
has been pursued in \cite{Miguelstability}, where the results
obtained here are applied to the study of the stability of
relative equilibria for simple mechanical systems at singular
momentum values. Another interesting future research direction is
to see how the results for the symplectic normal space obtained
here can advance towards a general explicit Symplectic Slice
Theorem for cotangent-lifted actions. In general the proof of this
theorem does not provide an explicit equivariant symplectomorphism
between the model space and the tubular neighborhood of the group
orbit in the Hamiltonian $G$-space. Only the existence of such a
map can be shown, together with its main properties. However,
recent results (see \cite{Schmah}) show that, at least for group
orbits consisting of points with totally isotropic momenta, such
an explicit construction is possible in the cotangent bundle case.
In \cite{Schmah}, using a completely different approach based on
singular commuting reduction, the expressions for the symplectic
normal space for cotangent-lifted actions at some particular types
of points are also computed, in agreement with our results in
subsections \ref{totallyisotropic} and \ref{relativeequilibria}.
Since we obtain a characterization of the symplectic normal space
at any point in $T^*Q$, it is expected that this will contribute
to generalizing this explicit symplectomorphism for
cotangent-lifted actions from totally isotropic to general
momentum values.
\paragraph{{\bf Notation:}} Unless otherwise specified, $Q$ will denote a
smooth and finite dimensional manifold. Throughout this paper we
use the symbol $\mathfrak{X}(Q)$ for the space of smooth vector
fields on $Q$ and
$\Lie_Y:\mathfrak{X}(Q)\rightarrow\mathfrak{X}(Q)$ for the Lie
derivative along $Y\in\mathfrak{X}(Q)$. If $i:V\hookrightarrow W$
is a linear injection of linear spaces, we denote by
$\Proj_V:W^*\rightarrow V^*$ its dual projection. Whenever a Lie
group $G$ acts smoothly on a space $X$ we denote by $G_x$ the
stabilizer of the element $x\in X$ and by $g\cdot x$ the action of
an element $g\in G$ on $x\in X$. We denote by $\g$ or
$\mathrm{Lie}\,(G)$ the Lie algebra of $G$ and by $\xi\cdot x$ or
$\xi_X(x)$ the infinitesimal action of $\xi\in\g$ on $X$. All
actions are assumed to be smooth. If $V$ is a linear space the
canonical pairing between $V^*$ and $V$ is denoted by
$\langle\cdot,\cdot\rangle$. Finally, the annihilator in $V^*$ of
a vector subspace $K\subset V$ is denoted by $K^\circ$.

\paragraph{{\bf Acknowledgements:}}
This research was mainly supported by the EU funding for the
Research Training Network MASIE, Contract No. HPRN-CT-2000-00113.
M. Perlmutter also thanks the Marsden Fund of the Royal Society of
New Zealand, as well as the Bernoulli Center of the EPFL
(Switzerland) for financial support. M.E. Sousa-Dias and M.
Perlmutter also thank the support the Mathematics Department of
IST and FCT (Portugal), in particular through the programs
POCTI/FEDER. We would like to thank Gianluigi del Magno for
pointing to us the importance of the Sasaki metric in the theory
of geodesic flows. Finally, we thank the anonymous referees whose
comments substantially improved the paper.
    \section{Hamiltonian actions and the symplectic normal space}\label{secWA}
 Let $(\PP,\w)$ be a symplectic manifold endowed with a
smooth and proper Hamiltonian action of a Lie group $G$ with Lie
algebra $\g$ and equivariant momentum map $\J:\PP\rightarrow
\mathfrak{g}^*$. Let $z\in\PP$ be a point with stabilizer
$G_z=\{g\in G\, :\, g\cdot z=z\}$ and momentum $\J (z)=\mu$ and
denote by $G_\mu$ the stabilizer of $\mu$ for the coadjoint
representation of $G$ and by $\g_\mu$ its Lie algebra. The
symplectic normal space at $z$ is the linear space $N=\ker
T_z\J/(\g_\mu\cdot z)$. The space $N$ is endowed with a natural
symplectic form defined as
$$\Omega (\left[u_1\right],\left[u_2\right]):=\w(z)(u_1,u_2),$$
where $u_1,u_2\in \ker T_z\J$. The induced linear $G_z$-action on
$T_z\PP$ descends to a linear Hamiltonian action on $N$, with
associated momentum map $\J_N:N\rightarrow \g_z^*$.

We obtain a realization of $N$ as a linear subspace of $T_z\PP$ in
the following way: Choose a $G_z$-invariant splitting
\be\label{introV}(\g\cdot z)^\w=\g_\mu\cdot z\oplus V,\ee where
$(\g\cdot z)^\w$ denotes the symplectic orthogonal of $\g\cdot z$.
The space $V$, equipped with the restricted symplectic form
$\w\rrestr{V}$, is a $G_z$-invariant maximal symplectic subspace
of $(\g\cdot z)^\w$. The group $G_z$ acts on $V$ in a Hamiltonian
fashion with associated equivariant momentum map
$\J_V:V\rightarrow \g_z^*$. Any such choice of invariant
complement $V$ provides a $G_z$-equivariant symplectomorphism
$N\simeq V$ that relates $\J_N$ and $\J_V$. Thus, all choices of
the subspace $V$ at $z$ are also symplectomorphic. We can think of
$N$ as the equivalence class of all such subspaces $V\subset
T_z\PP$. For this reason, hereafter we will also refer to both $V$
and $N$ as the symplectic normal space.

    \section{The isomorphism $T_{p_x}(T^*Q)\simeq
    T^*(T_xQ)$.}\label{secisom}

The setup for the rest of the paper is the following: Let
$\PP=T^*Q$ equipped with its canonical symplectic form, where
$(Q,\ll\cdot,\cdot\gg)$ is a Riemannian manifold. $G$ is a Lie
group acting on $Q$ properly by isometries and on $T^*Q$ by
cotangent lifts. This lifted action is also proper and
Hamiltonian, with equivariant momentum map
\be\label{cotangentmomentummap} \langle\J(p_x),\xi\rangle=\langle
p_x,\xi_Q(x)\rangle,\quad\forall p_x\in T_x^*Q\subset
T^*Q,\,\xi\in\g .\ee By construction, the bundle projection
$\tau:T^*Q\rightarrow Q$ is also equivariant. We are interested in
constructing the symplectic normal space $V$ at any point $p_x$.

 To achieve this objective, our first task is
to obtain convenient descriptions of the infinitesimal generators
of the cotangent-lifted action of $G$ on $T^*Q$ at $p_x$, and we
want these descriptions to incorporate the bundle structure of our
symplectic manifold $T^*Q$ into the description of the tangent
space $T_{p_x}(T^*Q)$. For that, we consider the Ehresmann
connection on $T^*Q$ associated to the Levi-Civita connection
$\nabla$ on $Q$, and the corresponding Whitney sum vector bundle
$$T(T^*Q)=\mathcal{H}(T^*Q)\oplus \mathcal{V}(T^*Q),$$ where
$\mathcal{H}(T^*Q)$ and $\mathcal{V}(T^*Q)$ are respectively the
horizontal and vertical bundles relative to this connection.

Let $p_x\in T^*Q$ be a given point over $x$. The connection map
$K:T_{p_x}(T^*Q)\rightarrow T_x^*Q$ is defined as follows: Let
$Y\in T_{p_x}(T^*Q)$ and let $\hat c(t)$ be a local curve in
$T^*Q$ such that $\hat c(0)=p_x$ and
$\frac{d\hat{c}(t)}{dt}\rrestr{t=0}=Y$. Then \be\label{K}
K(Y):=\frac{D^\nabla_c\hat{c}(t)}{Dt}\restr{t=0}\in T_x^*Q,\ee
i.e. the evaluation at time zero of the covariant differential
associated to $\nabla$ of the covector field $\hat c(t)$ along
$c(t)=\tau(\hat c(t))$. At $p_x$, the vertical space is given by
$\mathcal{V}_{p_x}:=\ker T_{p_x}\tau$, and the horizonal space by
$\mathcal{H}_{p_x}(T^*Q):=\ker K$. As the restrictions
$T_{p_x}\tau:\mathcal{H}_{p_x}(T^*Q)\rightarrow T_xQ$ and
$K:\mathcal{V}_{p_x}(T^*Q)\rightarrow T_x^*Q$ are isomorphisms, we
define a linear isomorphism $I:T_{p_x}(T^*Q)\rightarrow T_xQ\oplus
T_x^*Q$ by \be\label{Sasaki map} I(Y):=(T_{p_x}\tau (Y),K(Y)).\ee
We will call $I(Y)$ the $I$-representation of $Y$.

Since $G$ acts on $Q$ by isometries, the horizontal and vertical
distributions are $G$-invariant, and hence the map $I$ is
$G_{p_x}$-equivariant with respect to the induced linear
$G_{p_x}$-actions on $T_{p_x}(T^*Q),\,T_xQ$ and $T^*_xQ$. Another
way to see this is to realize $\mathcal{H}(T^*Q)$ as the
orthogonal complement to $\mathcal{V}(T^*Q)$ with respect to the
Sasaki metric on $T^*Q$, for which the lifted action is isometric
(see \cite{Sasaki}).

 We consider now finer properties of
the isomorphism $I$ due to the presence of the isometric
$G$-action. Let $p_x\in T_x^*Q$ and let $H=G_x$, with Lie algebra
$\h$. Since $G_{p_x}$ is compact by the properness of the lifted
$G$-action, we can choose a $G_{p_x}$-invariant splitting of $\g$
(see Section 6 for details) \be \label{ghr}\g=\h\oplus\lir .\ee
According to \eqref{ghr}, we can
 write  every element of $\g$ uniquely as $\eta=\eta^\h+\eta^\lir$,
with $\eta^\h\in\h$ and $\eta^\lir\in\lir$. Also, with respect to
$\ll\cdot,\cdot\gg$ the tangent space $T_xQ$ is split orthogonally
as $T_xQ=\g\cdot x\oplus \Sl$, where $\Sl=(\g\cdot x)^\perp$ is a
linear slice  at $x$ for the $G$-action on $Q$.

In this paper we will consider the following linear actions of
several subgroups of $G$: The group $G_{p_x}$ acts on $\lir$ and
$\lir^*$ by restricting the adjoint and coadjoint representations
of $G$. The group $H$ acts on $\Sl$ by its induced representation
on $T_xQ$ and on $\Sl^*$ by its contragredient representation,
$\langle h\cdot\beta,b\rangle=\langle\beta,h^{-1}\cdot b\rangle$,
for $h\in H,\,b\in\Sl$ and $\beta\in\Sl^*$. The infinitesimal
generator map $\g\rightarrow \g\cdot x$, $\xi\mapsto \xi_Q(x)$ is
not in general an isomorphism since it may have nontrivial kernel
$\h$. However, by~\eqref{ghr} this map is an isomorphism if its
domain is restricted to $\lir$. The locked inertia tensor,
introduced next, is a useful family of bilinear forms on $\g$ that
will help us to express the $I$-representation of vectors in
$T_{p_x}(T^*Q)$. It can be seen as a family of (degenerate)
metrics on $\g$ induced from the Riemannian metric on $Q$ by the
infinitesimal generator map.
\begin{defn}\label{deflock} The locked inertia tensor on a Riemannian manifold
$(Q,\ll \cdot,\cdot\gg)$ equipped with an isometric action of a
Lie group $G$ is the map that associates to each point $x\in Q$
the symmetric bilinear form on $\mathfrak{g}$ defined by \be \II
(x)(\xi,\eta):=\ll \xi_Q(x),\eta_Q(x)\gg . \ee
\end{defn}
The following well known property of the locked inertia tensor
will be used throughout this paper. Its proof can be found in
chapter 5 of \cite{MarLec}.

\begin{lemma}\label{lockinertia}
For every $\xi,\eta,\lambda\in\g$, $g\in G$ and $x\in Q$ the
following identities hold. $$\begin{array}{l}\II(g\cdot
x)(\Ad_g\eta,\Ad_g\lambda)  = \II(x)(\eta,\lambda),\quad\mathrm{and}\vspace{3mm}\\
(\D\II\cdot\xi_Q(x))(\eta,\lambda)+\II
(x)(\ad_\xi\eta,\lambda)+\II
(x)(\eta,\ad_\xi\lambda)=0,\end{array}$$ where $\Ad$ and $\ad$
denote respectively the adjoint representations of $G$ and $\g$.
\end{lemma}

 The locked inertia tensor $\II (x)$, seen as a
bilinear symmetric two-form on $\g$, is degenerate with kernel
equal to the Lie algebra $\h=\g_x$. Consequently, its restriction
$$\rII :=\II (x)\rrestr{\lir\times\lir}:\lir\times\lir
\rightarrow\R$$ is a well-defined inner product on $\lir$. In
order to economize notation, we will denote by the same symbols
the maps $\II (x):\g\rightarrow\g^*$ and
$\rII:\lir\rightarrow\lir^*$ induced respectively by $\II (x)$ and
$\rII$. Since $T_xQ=\g\cdot x\oplus\Sl$ and as $\g\cdot x$ is
isomorphic to $\lir$ through the infinitesimal generator map
$\xi\mapsto \xi_Q(x)$, we can form the $G_{p_x}$-equivariant
linear isomorphisms \be\label{splittangent}
  T_xQ \simeq \lir\oplus \Sl \quad \mathrm{and}\quad
  T_x^*Q  \simeq  \lir^* \oplus \Sl^*, \ee
with the first isomorphism
  $\lir\oplus \Sl\rightarrow T_xQ$ defined by
  $(\xi,a)\mapsto \xi_Q(x)+a,$
and the second, $T^*_xQ\simeq \lir^*\oplus \Sl^*$ obtained by
dualization.

In this representation, the metric Legendre map $\FL
:T_xQ\rightarrow T_x^*Q$, which associates to each vector $v_x\in
T_xQ$ the covector $\ll v_x,\cdot\gg\in T_x^*Q$, is written as
\be\label{FLsasaki} \left<\FL(\xi_1,a_1),(\xi_2,a_2)\right> = \II
(x) (\xi_1,\xi_2)+\ll a_1,a_2\gg_\Sl =  \rII (\xi_1,\xi_2)+\ll
a_1,a_2\gg_\Sl ,\ee where $\xi_1,\xi_2\in\lir$, $a_1,a_2\in\Sl$
and  $\ll\cdot,\cdot\gg_\Sl$ denotes the restriction to $\Sl$ of
the inner product in $T_xQ$ induced by the metric. Note that
$\lir\simeq \g/\h$ and $\lir^*\simeq \h^\circ$.

Finally, the composition of the two isomorphisms
~\eqref{splittangent} with  $I$ defined in \eqref{Sasaki map}
allows us to identify $G_{p_x}$-isomorphically
$T_{p_x}(T^*Q)\simeq (\lir\oplus \Sl )\oplus (\lir^*\oplus
\Sl^*)$. Therefore, the $I$-representation of a vector $v_{p_x}$
tangent to $T^*Q$ at $p_x$ is given by the quadruple
$$I(v_{p_x})=(\xi,a;\nu,\alpha)\in (\lir\oplus \Sl) \oplus (\lir^*\oplus
\Sl^*),$$ where $\xi,a,\nu$ and $\alpha$ are uniquely defined by
\begin{eqnarray*}T_{p_x}\tau (v_{p_x}) & = & \xi_Q(x)+a\in T_xQ,\quad\mathrm{and}\\ K (v_{p_x}) &
=& \FL
\left(\left(\rII^{-1}(\nu)\right)_Q(x)\right)+\alpha.\end{eqnarray*}

\paragraph{{\bf Remarks:}}
\begin{enumerate}\item There is no loss of generality in
supposing that $Q$ is a Riemannian manifold and that $G$ acts
isometrically on it as long as $Q$ is paracompact since, by
properness of the $G$-action, one can always find a Riemannian
structure on $Q$ invariant under the given action (see
\cite{DuiKol}). \item A trivial but extremely important
observation when working with cotangent-lifted actions is that,
since both $\tau:T^*Q\rightarrow Q$ and $\J:T^*Q\rightarrow \g^*$
are $G$-equivariant, for a point $p_x\in T^*Q$ with base point
$x=\tau (p_x)$ and momentum $\mu=\J (p_x)$, the following relation
among the three involved isotropy groups holds: $G_{p_x}\subset
G_x\cap G_\mu$.
\end{enumerate}

The canonical symplectic form of $T^*Q$ has a particularly simple
expression in the representation provided by the map $I$.
\begin{lemma}
Let $\w$ denote the canonical symplectic form on $T^*Q$ and $Y_1,
Y_2\in T_{p_x}(T^*Q)$. Then
\be\label{sasakisympl1}\w(p_x)(Y_1,Y_2)=\langle K(Y_2),T_{p_x}\tau
(Y_1)\rangle-\langle K(Y_1),T_{p_x}\tau (Y_2)\rangle .\ee
\end{lemma}
\begin{proof}
Let $n=\dim Q$ and let $(x^1,\ldots,x^n,p_1,\ldots,p_n)$ be local
coordinates of a bundle chart $U\times \R^n$ of $T^*Q$ with
$U\subset Q$. Then local frames for $\mathcal{V}(T^*U)$  and
$\mathcal{H}(T^*Q)$ are given respectively by
\begin{eqnarray*}\frac{\partial}{\partial p_i} & & \quad\mathrm{and}\\
\frac{\delta}{\delta x^i} & = & \frac{\partial}{\partial
x^i}+\Gamma_{ij}^kp_k\frac{\partial}{\partial p_j},\quad i,j,k\in
(1,\ldots,n).\end{eqnarray*} We have used the Einstein convention
for index summation. Here, the $\Gamma$'s are the Christoffel
symbols of $\nabla$, defined by $\nabla_\frac{\partial}{\partial
x^i}\frac{\partial}{\partial
x^j}=\Gamma_{ij}^k\frac{\partial}{\partial x^k}$ or
$\nabla_\frac{\partial}{\partial x^i}dx^j=-\Gamma_{ik}^j dx^k$.
 Then, if $Y\in
T_{p_x}(T^*Q)$ is written as $Y=A^i\frac{\delta}{\delta
x^i}+B_i\frac{\partial}{\partial p_i}$, $T_{p_x}\tau
(Y)=A^i\frac{\partial}{\partial x^i}$ and $K(Y)=B_idp^i$. Recall
that in these local coordinates, $\w$ is characterized by $
\w(\frac{\partial}{\partial x^i},\frac{\partial}{\partial
p_j})=\delta_i^j$ and $\w(\frac{\partial}{\partial
x^i},\frac{\partial}{\partial x^j})=\w(\frac{\partial}{\partial
p_i},\frac{\partial}{\partial p_j})=0$. Therefore we have
$\w(\frac{\delta}{\delta x^i},\frac{\partial}{\partial
p_j})=\delta_i^j$. Also $\w(\frac{\delta}{\delta
x^i},\frac{\delta}{\delta
x^j})=p_k(\Gamma_{ji}^k-\Gamma_{ij}^k)=0$ since $\nabla$ has zero
torsion. Now \eqref{sasakisympl1} follows by bilinearity of $\w$.
\end{proof}
 Note that the resulting symplectic form $\w$ is still given by \eqref{sasakisympl1} if we perform the
same construction with respect to the Ehresmann connection
associated to any affine connection on $Q$ with zero torsion.

As a consequence of~\eqref{sasakisympl1}, the symplectic form at
$p_x$ in the four-way $I$-representation of $T_{p_x}(T^*Q)$ is
expressed as \be\label{sasakisymplectic}
\w(p_x)\left((\xi_1,a_1;\nu_1,\alpha_1),(\xi_2,a_2;\nu_2,\alpha_2)\right)
=\left<\nu_2,\xi_1\right>-\left<\nu_1,\xi_2\right>+
\left<\alpha_2,a_1\right>-\left<\alpha_1,a_2\right> .\ee

 \section{Local vector fields on $Q$}\label{secvectorfields}
In order to compute the map $K$ applied to several types of
vectors belonging to $T_{p_x}(T^*Q)$, we need to use a family of
locally defined vector fields on $Q$ such that their restrictions
at $x$ span $T_xQ$.  Such a  family should be adapted in some
sense to the $G$-action. This section is devoted to the
construction of such vector fields and to the development of a
series of results about them, summarized in Lemma \ref{orto},
Proposition \ref{brackets} and Proposition \ref{covariant}.

We will work in a local model of $Q$ given by a neighborhood $O$
of $[e,0]$ in the associated bundle $G\times_H\Sl$ for which the
map $\phi:G\times_H\Sl\rightarrow Q$ given by
$\phi([g,s])=g\cdot\exp_x s$ restricts to a local $G$-equivariant
diffeomorphism onto a $G$-invariant neighborhood $U$ of $G\cdot x$
according to Palais' Tube Theorem (see
\cite{DuiKol,Koszul,OrRa2003,Palais}). Here $\exp_x$ is the
Riemannian exponential at $x$ and $G\times_H\Sl$ is the orbit
space of $G\times\Sl$ by the twisted $H$-action
$h\cdot(g,s)=(gh^{-1},h\cdot s)$, and $G$ acts on $G\times_H\Sl$
by $g'\cdot[g,s]=[g'g, s]$.

Note that,  as  the map $\pi_H :G\times \Sl\rightarrow G\times_H
\Sl$ defines a principal bundle,  we can write $O=\overline{O}/H$,
where $\overline{O}\subset G\times\Sl$ is an $H$-saturated
neighborhood of $(e,0)$. We now define an $H$-equivariant flow on
$\overline{O}$. Consider any inner product on $\g$
 such that $\lir$ in the splitting \eqref{ghr} is obtained as
 $\lir=\h^\perp$. This inner
product can be extended  to a right-invariant Riemannian metric on
$G$, and hence invariant for the $H$-action on $G$ given by
$h\cdot g= gh^{-1}$, for every $h\in H$ and $g\in G$. Note that
the tangent space to $G$ at the identity is  $\g$ and  the
$H$-orbit through $e$ is exactly $H$, with tangent space $\h$.
Therefore $\lir$ is an orthogonal linear slice for the $H$-action
at $e$ with respect to this metric. Denote by $\exp_e
:T_eG=\g\rightarrow G$ the associated Riemannian exponential at
the identity (not to be confused with the group exponential
$\Exp$). A straightforward application of the Tube Theorem
guarantees that there exists an $H$-invariant neighborhood of $H$
in $G$ in which every element $g$ can be written uniquely as
$$g=(\exp_e \xi^\lir)h^{-1},$$ for some pair $(h,\xi^\lir)\in
H\times\lir$. Uniqueness follows from the fact that the action of
$H$ on $G$ is free, and hence the associated bundle playing the
role of the tube is just the direct product $H\times\lir$.

\begin{lemma}\label{lvectorfields}  Let $\overline{O}\subset G\times \Sl$ be
a small enough $H$-invariant neighborhood of $(e,0)$ for which any
$(g,s)\in \overline{O}$ is written uniquely as
$((\exp_e\xi^\lir)h^{-1},s)$ for some $\xi^\lir\in\lir$ and $h\in
H$. Then, for every $v\in\Sl$ and $t$ small enough,
 $\overline{F}_v^t(g,s):=(g,s+th\cdot v)$ defines a local
$H$-equivariant flow on $\overline{O}$.
\end{lemma}
\begin{proof}
We need to check that $\overline{F}_v^{t}$ is a well defined
one-parameter group of local diffeomorphisms. Obviously
$\overline{F}_v^0=\mathrm{id}_{\overline{O}}$, and
$$\begin{array}{lll}
(\overline{F}_v^{t_1}\circ\overline{F}_v^{t_2})(g,s) & = &
\overline{F}_v^{t_1}(g,s+t_2h\cdot v)= (g,s+t_2h\cdot v+t_1h\cdot
v)\\& = & (g,s+(t_1+t_2)h\cdot v)=\overline{F}_v^{t_1+t_2}(g,s).
\end{array}$$
Furthermore, $\overline{F}_v^{t}$ is $H$-equivariant, since for
every $h'\in H$
$$\begin{array}{lll}
\overline{F}_v^t (h'\cdot (g,s)) & = & \overline{F}_v^t
(gh'^{-1},h'\cdot s)= \overline{F}_v^t((\exp_e
t\xi^\lir)h^{-1}h'^{-1},h'\cdot s)\\ & = & \left((\exp_e
t\xi^\lir)(h'h)^{-1},h'\cdot s+t(h'h)\cdot v\right)\\  & = &
h'\cdot ((\exp_e t\xi^\lir)h^{-1},s+th\cdot v)=
h'\cdot(\overline{F}_v^t (g,s)).
\end{array}$$\end{proof}
If necessary we can shrink $\overline{O}$ in the statement of
Lemma \ref{lvectorfields} in order to have $O=\overline{O}/H$
inside the domain of injectivity of the tube map
$\phi:G\times_H\Sl\rightarrow U$. Since $\overline{F}_v^t$ defines
an $H$-equivariant flow on $\overline{O}$, it descends to a flow
on $O\subset G\times_H\Sl$, and so applying $\phi$ to it we have a
well defined flow $F_v^t$ on $U\subset Q$. It has an associated
vector field $\overline{v}$ obtained by differentiating $F_v^{t}$.
Therefore we have
\begin{prop}\label{vector field definitions} Let $\overline{O}$ be as in
Lemma~\ref{lvectorfields}, and let $U\subset Q$ be the image under
$\phi:G\times_H\Sl\rightarrow Q$ of $O=\pi_H(\overline{O})$.
Shrink $\overline{O}$ if necessary so that
$\phi\rrestr{O}:O\rightarrow U$ is injective. Then every $x'\in U$
can be uniquely written as $x' = \phi ( [g,s])= \phi
([(\exp_e\xi^\lir)h^{-1},s])$ for $\xi^\lir\in\lir$ and $h\in H$.
Furthermore, for every $v\in\Sl$ and $\eta\in\g$ the following
formulae define vector fields $\overline{eta}$ and $\overline{v}$
on $U$: \be\label{campo1}
  \overline{\eta}(x') := \eta_Q(x'),
\ee \be\label{campo2}\begin{array}{ll}
 \overline{v}(x')  := \frac{d}{dt}\rrestr{t=0}
F_v^t(x')&= \frac{d}{dt}\rrestr{t=0}\phi ([g,s+th\cdot v])\vspace{1mm} \\
&= \frac{d}{dt}\rrestr{t=0}\phi ([(\exp_e\xi^\lir)h^{-1} ,
s+th\cdot v]).
\end{array}\ee
\end{prop}
\paragraph{{\bf Remark:}} Note that  the vector fields $\overline{\eta}$ and
$\overline{v}$   at the point $x= \phi([e,0])$ satisfy
$\overline{\eta}(x) = \xi_Q(x)$ and $\overline{v}(x) =  v$, for
$\eta\in\g$ and $v\in \Sl$. As $\lir\simeq \{\eta_Q(x)\,:\,\forall
\eta\in\g\}$, then the family of vector fields in Proposition
~\ref{vector field definitions} evaluated at $x$ spans
 $T_xQ$.

Now we establish some properties of the vector fields defined by
\eqref{campo1} and \eqref{campo2}. The next result just gives
another way of computing $\overline{\eta}(x')=\eta_Q(x')$ at a
point near $x$ for an element $\eta\in\g$. This alternative
construction will be used later.
\begin{lemma}\label{alternativo}
Let $U$ be as in the Proposition~\ref{vector field definitions}, $
x'=\phi([g,s])\in U$ for $g\in G$ and $s\in \Sl$, and $\eta\in\g$.
Then \be \overline{\eta}(x')=\frac {d}{dt}\restr{t=0}\phi([(\exp_e
t\eta)g,s]). \ee
\end{lemma}
\begin{proof}
By \eqref{campo1}, $\eta_Q(x')$ is just the infinitesimal
generator for the $G$ action on $Q$ corresponding to $\eta$, that
is $\eta_Q(x')  =  \frac{d}{dt}\rrestr{t=0}(\Exp t\eta)\cdot x'$.
As the map $\phi$ is $G$-equivariant, then for $\eta\in\g$ we have
$$\begin{array}{lll}
\eta_Q(x') & = & \frac{d}{dt}\rrestr{t=0}(\Exp t\eta)\cdot\phi
([g,s])=\frac{d}{dt}\rrestr{t=0}\phi ((\Exp t\eta)\cdot[g,s])\vspace{1mm}\\
 & = & \frac{d}{dt}\rrestr{t=0}\phi ([(\exp_e t\eta) g,s]),
\end{array}$$
 where the last equality follows from the fact that both curves, $g(t)=\exp_e
t\eta$ and $ g(t)=\Exp t\eta$, are tangent to $\eta$ at $t=0$.
\end{proof}
The next lemma shows that the vector fields~\eqref{campo1} and
\eqref{campo2} are orthogonal along $G\cdot x$.
\begin{lemma}\label{orto} Let $U=\phi (O)$ be as in Proposition~\ref{vector
field definitions}. Then, for all $x'\in G\cdot x\cap U,\,v\in\Sl$
and $\xi\in\g$ we have $\ll\overline{v}(x'),\overline{\xi}(x')\gg =0$.
\end{lemma}
\begin{proof}Let $x'=g\cdot x=\phi([g,0])$ and $g=(\exp_e \xi^\lir
)h^{-1}$ for $h\in H=G_x$. Then
 $$\overline{v}(x')=\frac {d}{dt}\restr{t=0} \phi
([g,th\cdot v])=g\cdot\left(\frac {d}{dt}\restr{t=0}\phi
([e,th\cdot v])\right)=g\cdot (h\cdot v)=(\exp_e\xi^\lir)\cdot
v.$$ Since $\ll\cdot,\cdot\gg$ is $G$-invariant,
$(\exp_e\xi^\lir)\cdot v\in(\g\cdot ((\exp_e\xi^\lir)\cdot
x))^\perp=(\g\cdot x')^\perp$.
\end{proof}
The next proposition on the evaluation at $x$ of   the Lie
brackets of the previously defined vector fields is of fundamental
importance for the study of the Levi-Civita connection of $Q$.
\begin{prop}\label{brackets}
The commutators at $x$ of the local vector fields \eqref{campo1}
and \eqref{campo2}  satisfy
\begin{itemize} \item[a)]
$\left[\overline{\xi}_i,\overline{\xi}_j\right](x)=-\overline{\left[\xi_i,\xi_j\right]}(x)$,
\item[b)] $\left[\overline{v_a},\overline{v_b}\right](x)=0$,
\item[c)]
$\left[\overline{\lambda},\overline{v}\right](x)=\left\{\begin{array}{ccc}
  0 & \mathrm{if} & \lambda\in\lir \\
  -\lambda\cdot v & \mathrm{if} & \lambda\in\h ,
\end{array}\right.$
\end{itemize}
for every $\lambda,\,\xi_i,\,\xi_j\in\g$ and
$v,\,v_a,\,v_b\in\Sl$.
\end{prop}
\begin{proof}
Item a) is just the rephrasing of the well known identity
$[\xi_Q,\eta_Q]=-([\xi,\eta])_Q$ for left actions.

For b), the result is a consequence of the computation of the Lie
derivative as follows.
$$\begin{array}{ll} \left[\overline{v_a},\overline{v_b}\right](x)& =
\Lie_{\overline{v_a}}\overline{v_b}(x)  =
\frac{d}{dt}\rrestr{t=0}\left((F^t_{v_a})^*\overline{v_b}\right)(x)\vspace{1mm}\\&
=
\frac{d}{dt}\rrestr{t=0}T_{F^{t}_{v_a}(x)}F^{-t}_{v_a}\left(\overline{v_b}
(F^t_{v_a}(x))  \right)\vspace{1mm}\\

& =  \frac{d}{dt}\rrestr{t=0}\frac{d}{ds}\rrestr{s=0}
F_{v_a}^{-t}\left(F_{v_b}^s(F_{v_a}^t(x))\right)\vspace{1mm}\\ & =
\frac{d}{dt}\rrestr{t=0}\frac{d}{ds}\rrestr{s=0}\phi
([e,-tv_a+sv_b+tv_a])=0.
\end{array}$$
For c), let us first consider $\lambda\in\lir$, then
$$\begin{array}{ll} -[\overline{\lambda},\overline{v}](x) &  =
\Lie_{\overline{v}}\overline{\lambda}(x)  =
\frac{d}{dt}\rrestr{t=0}\left((F^t_{v})^*\overline{\lambda}\right)(x)\\
& =
\frac{d}{dt}\rrestr{t=0}T_{F^{t}_{v}(x)}F^{-t}_{v}\left(\overline{\lambda}
(F^t_{v}(x))  \right)\vspace{1mm}\\
& = \frac{d}{dt}\rrestr{t=0}\frac{d}{ds}\rrestr{s=0}
F_v^{-t}\left(\phi([\exp_e s\lambda,tv])\right)\vspace{1mm}\\ & =
\frac{d}{dt}\rrestr{t=0}\frac{d}{ds}\rrestr{s=0}\phi([\exp_e
s\lambda,-tv+tv])=0,
\end{array}$$
where the fourth equality follows from Lemma \ref{alternativo}.\\
Analogously, if $\lambda\in\h$ we have
$$\begin{array}{ll} [\overline{v},\overline{\lambda}](x)
&=\Lie_{\overline{v}}\,\overline{\lambda}(x)  =
\frac{d}{dt}\rrestr{t=0}\frac{d}{ds}\rrestr{s=0}
F_v^{-t}\left(\phi([\Exp (s\lambda),tv])\right)\vspace{1mm}\\ & =
\frac{d}{dt}\rrestr{t=0}\frac{d}{ds}\rrestr{s=0} \phi ([\Exp
(s\lambda),tv-t(\Exp (s\lambda))^{-1}\cdot v)])\vspace{1mm}\\  & =
\frac{d}{dt}\rrestr{t=0}\frac{d}{ds}\rrestr{s=0} \Exp (s\lambda)
\cdot \exp_x(tv-t(\Exp (-s\lambda))\cdot v)).
\end{array}$$
Now we compute the above expression, which is of the general form
$$\frac{d}{dt}\restr{t=0}\frac{d}{ds}\restr{s=0}g(s)\cdot\exp_x tk(s),$$
where $g(s)=\Exp (s\lambda)$ is a curve in $G$ and $k(s)=v-\Exp(-
s\lambda)\cdot v$ is a curve
 in $\Sl$. We can then write
$$\frac{d}{dt}\restr{t=0}g(s)\cdot\exp_x tk(s)=g(s)\cdot T_0\exp_x
k(s)=g(s)\cdot k(s),$$ since $T_0\exp_x=\mathrm{id}_{T_xQ}$.
Finally, from the expressions of $g(s)$ and $k(s)$ we find that
$$g(s)\cdot k(s)=\Exp (s\lambda)\cdot v-v,$$
and then finally $
[\overline{v},\overline{\lambda}](x)=\frac{d}{ds}\restr{s=0}g(s)\cdot
k(s)=\lambda\cdot v.$
\end{proof}

The next result will play a key role in obtaining the
$I$-representation of vectors tangent to $T^*Q$, simplifying
considerably the computation of the connection map $K$.

\begin{prop}\label{covariant} For any element $\eta\in\g$ consider the
unique decomposition $\eta=\eta^\h+\eta^\lir$ with respect to the
splitting \eqref{ghr}. Let $v,w\in\Sl$,
$\lambda,\xi,\xi_i,\xi_j\in\g$. Then

\begin{enumerate}
\item $ \ll
\nabla_{\overline{\xi_i}}\overline{\xi_j}(x),\lambda_Q(x)\gg  =
\frac 12\left\{ \left(\mathbf{D}\II \cdot
(\xi^\lir_i)_Q(x)\right)(\xi_j,\lambda)-\II
(x)(\xi^\lir_i,\left[\xi_j,\lambda\right])\right\}$.

\item $\ll
\nabla_{\overline{\xi_i}}\overline{\xi_j}(x),w\gg=-\frac
12\left(\mathbf{D}\II\cdot w\right)(\xi_i^\lir,\,\xi_j^\lir)$.

\item $\ll
\nabla_{\overline{\xi}}\,\overline{v}(x),\lambda_Q(x)\gg=\frac
12\left(\mathbf{D}\II\cdot v\right)(\xi^\lir ,\,\lambda)$.

\item $\ll
\nabla_{\overline{v}}\,\overline{\xi}(x),\lambda_Q(x)\gg=\frac
12\left(\mathbf{D}\II \cdot v\right)(\xi^\lir,\lambda)$.

\item $\ll
\nabla_{\overline{v}}\,\overline{\xi}(x),w\gg=\ll\nabla_{\overline{\xi}}\,\overline{v}(x),w\gg+\ll\xi^\h\cdot
v,w\gg_\Sl $.

\end{enumerate}
\end{prop}

\paragraph{{\bf Notation:}} We will introduce a concise notation for
$\nabla_{\overline{\xi}}\,\overline{v}(x)$. Recall that the map
$(\bar{\xi},\bar{v})\mapsto
\nabla_{\overline{\xi}}\,\overline{v}(x)$ is linear in the
argument $\overline{\xi}$, and depends on $\overline{\xi}$ only
through its value at $x$, so we can write
\be\label{defC}\ll\nabla_{\overline{\xi}}\,\overline{v}(x),w\gg=\ll
C(v)(\xi^\lir ),w\gg_\Sl .\ee This defines the bilinear map $\ll
C(v)(\cdot),\cdot\gg_\Sl:\lir\times\Sl \rightarrow\R$, where
$\overline{v}$ is determined by $v$ through \eqref{campo2}.

\begin{proof} We will choose extensions
$\overline{\lambda}$ and $\overline{w}$ of $\lambda_Q(x)$ and
$w\in\Sl$, respectively. Since, for every $X,Y\in\mathfrak{X}(Q)$,
the assignment $(X,Y)\mapsto \nabla_XY$ is $C^\infty(Q)$-linear in
$X$, then $\nabla_{\overline{\xi}}Y(x)=0$ for every $\xi\in\h$,
and therefore, the $\h$-component of elements of $\g$ does not
appear in items (1) to (4). Also, since $\nabla$ has zero torsion,
(2) is symmetric under the permutation $\xi_i\leftrightarrow\xi_j$
and hence it does not depend on $\xi_j^\h$.

 For (1), note that $\xi_j=\xi_j^\lir+\xi_j^\h$ and that $\nabla_XY$ is
linear in $Y$, and therefore
$$\nabla_{\overline{\xi_i}}\overline{\xi_j} (x)=\nabla_{\overline{\xi_i^\lir}}\overline{\xi_j}(x)+
\nabla_{\overline{\xi_i^\h}}\overline{\xi_j}(x)=
\nabla_{\overline{\xi^\lir_i}}\overline{\xi_j}(x)=
\nabla_{\overline{\xi^\lir_i}}\overline{\xi^\lir_j}(x)+\nabla_{\overline{\xi^\lir_i}}\overline{\xi^\h_j}(x)
.$$ Let us compute
$\nabla_{\overline{\xi^\lir_i}}\overline{\xi^\h_j}$ and
$\nabla_{\overline{\xi^\lir_i}}\overline{\xi^\lir_j}$ separately.
For the $\lir -\h$ part, using the standard formula for the
covariant derivative of the Levi-Civita connection, we have
$$\begin{array}{lll}
2\ll\nabla_{\overline{\xi^\lir_i}}\overline{\xi^\h_j}(x),\lambda_Q(x)\gg
& = &
\overline{\xi^\lir_i}\left(\ll\overline{\xi^\h_j},\lambda_Q\gg\right)
(x)+\overline{\xi^\h_j}\left(\ll\overline{\xi^\lir_i},\lambda_Q\gg\right)
(x)\\& & -\overline\lambda \left(\ll
\overline{\xi^\lir_i},\overline{\xi^\h_j}\gg\right) (x) + \ll
\lambda_Q(x),[\overline{\xi^\lir_i},\overline{\xi^\h_j}](x)\gg\\ &
& +\ll
\overline{\xi^\h_j}(x),[\overline{\lambda},\overline{\xi^\lir_i}](x)\gg
 -\ll
\overline{\xi^\lir_i} (x),[\overline{\xi^\h_j},\overline{\lambda}
](x)\gg,
\end{array}$$
where the second and fifth terms vanish since
$\overline{\xi^\h_j}(x)=0$. In view of Definition \ref{deflock},
$$\begin{array}{lll} 2\ll
\nabla_{\overline{\xi^\lir_i}}\overline{\xi^\h_j}(x),\lambda_Q(x)\gg
& = & \left(\D\II\cdot
(\xi_i^\lir)_Q(x)\right)(\xi_j^\h,\lambda)-(\D\II\cdot\lambda_Q(x))(\xi_i^\lir,\xi_j^\h)\\
 & \quad & + \II (x)(\lambda,\ad_{\xi_j^\h}\,\xi_i^\lir)-\II
(x)(\xi_i^\lir,\ad_\lambda\xi_j^\h),\end{array}$$ where we have
applied Proposition~\ref{brackets}-a). Using Lemma
\ref{lockinertia} we obtain,
$$(\D\II\cdot\lambda_Q(x))(\xi_i^\lir,\xi_j^\h)+\II
(x)(\xi_i^\lir,\ad_\lambda\xi_j^\h)=-\II(x)(\ad_\lambda\xi_i^\lir,\xi_j^\h)
=0,\quad\mathrm{so}$$
$$\begin{array}{lll} 2\ll
\nabla_{\overline{\xi^\lir_i}}\overline{\xi^\h_j}(x),\lambda_Q(x)\gg
& = & (\D\II\cdot (\xi_i^\lir)_Q(x))(\xi_j^\h,\lambda)+\II
(x)(\lambda,\ad_{\xi_j^\h}\xi_i^\lir)\\
& = & (\D\II\cdot (\xi_i^\lir)_Q(x))(\xi_j^\h,\lambda)-\II
(x)(\xi^\lir_i,\ad_{\xi_j^h}\lambda ),
\end{array}$$
where the last equality also follows from Lemma \ref{lockinertia}.

 Now, with the same reasoning as for the $\lir -\lir$ part,
$$\begin{array}{lll} 2\ll
\nabla_{\overline{\xi^\lir_i}}\overline{\xi^\lir_j}(x),\lambda_Q(x)\gg
& = & (\D\II\cdot (\xi_i^\lir)_Q(x))(\xi_j^\lir,\lambda)  +
(\D\II\cdot
(\xi_j^\lir)_Q(x))(\xi_i^\lir,\lambda)\\
& \quad & -(\D\II\cdot \lambda_Q(x))(\xi_i^\lir,\xi_j^\lir)-\II
(x)(\lambda,\ad_{\xi_i^\lir}\xi_j^\lir)\\
& \quad  & -\II (x)(\xi_j^\lir,\ad_\lambda\xi_i^\lir)+\II
(x)(\xi_i^\lir,\ad_{\xi_j^\lir}\lambda)\\
& = & (\D\II\cdot (\xi_i^\lir)_Q(x))(\xi_j^\lir,\lambda)-\II
(x)(\xi_i^\lir,\ad_{\xi_j^\lir}\lambda) .\end{array}$$ We get (1)
by adding the $\lir -\h$ and the $\lir -\lir$ parts. For (2) note
that
$$\begin{array}{lll}
2\ll\nabla_{\overline{\xi_i}}\overline{\xi_j}(x),w\gg & = &
\overline{\xi^\lir_i}(\ll
\overline{\xi_j^\lir},\overline{w}\gg)(x)+\overline{\xi^\lir_j}(\ll
\overline{\xi^\lir_i},\overline{w}\gg)(x)\\
& \quad &
-\overline{w}(\ll\overline{\xi^\lir_i},\overline{\xi^\lir_j}\gg
)(x)+\ll w,[\overline{\xi^\lir_i},\overline{\xi^\lir_j}](x)\gg\\
& \quad & +\ll
(\xi_j)_Q(x),[\overline{w},\overline{\xi_i^\lir}](x)\gg -\ll
(\xi^\lir_i)_Q(x),[\overline{\xi^\lir_j},\overline{w}](x)\gg .
\end{array}$$
In the above expression the first two terms vanish by Lemma
\ref{orto} while, for the third, the use of Definition
\ref{deflock} yields
$$\overline{w}(\ll\overline{\xi^\lir_i},\overline{\xi^\lir_j}\gg
)(x)=(\D\II\cdot w)(\xi^\lir_i,\xi^\lir_j).$$ Finally, using
Proposition \ref{brackets}, and recalling that $w\in\Sl=(\g\cdot
x)^\perp$, one checks that the last three terms also vanish. To
prove (3), we expand its expression as
$$\begin{array}{lll}2\ll\nabla_{\overline{\xi}}\overline{v}(x),\lambda_Q(x)\gg &
= & \overline{\xi^\lir}(\ll\overline
v,\overline\lambda\gg)(x)+\overline
v(\ll\overline{\xi^\lir},\overline\lambda\gg)(x)\\
& & -\overline\lambda(\ll\overline{\xi^\lir},\overline
v\gg)(x)+\ll\lambda_Q(x),[\overline{\xi^\lir},\overline{v}](x)\gg\\
& & +\ll
v,[\overline\lambda,\overline{\xi^\lir}](x)\gg-\ll\xi^\lir_Q(x),[\overline
v ,\overline\lambda](x)\gg .\end{array}$$ The first and third
contributions vanish by Lemma~\ref{orto}, the fourth and sixth
vanish by Proposition~\ref{brackets}, and the fifth vanishes since
the commutator of infinitesimal generators is also an
infinitesimal generator, and hence perpendicular to $v$. Then the
only non-vanishing term is the second, and (3) is proved.

For (4), recall that $\nabla$ has zero torsion and so one can
write, for every $a\in T_xQ$
$$\ll\nabla_{\overline{v}}\,\overline{\xi}(x),a\gg=\ll\nabla_{\overline{\xi}}\,\overline{v}
(x),a\gg+ \ll [\overline{v},\overline{\xi}](x),a\gg.$$ By
Proposition \ref{brackets}  $\ll
[\overline{v},\overline{\xi}](x),\lambda_Q(x)\gg=0$. Consequently,
$$\ll\nabla_{\overline{v}}\,\overline{\xi}(x),\lambda_Q(x)\gg=\ll\nabla_{\overline{\xi}}\,\overline{v}
(x),\lambda_Q(x)\gg=\frac 12 \left(\D\II\cdot v\right)
(\xi^\lir,\lambda),$$ where the last equality follows by (3).

For (5), use again the expression for the torsion of $\nabla$ and
note that by Proposition ~\ref{brackets}
 $$\ll
[\overline{v},\overline{\xi}](x),w\gg=\ll\xi^\h\cdot v,w\gg.$$
\end{proof}

We end this section with a technical result involving the locked
inertia tensor and the map $C$. This will be needed in Section 6.

\begin{lemma}\label{intermch3}
Let $h\in G_{p_x},\,\xi\in\g$, $\eta\in\lir$ and $v,w\in T_xQ$,
then
\begin{enumerate}\item $ \Ad_{h^{-1}}^*\left[(\D\II\cdot v)(\xi)\right] =
(\D\II\cdot (h\cdot v))(\Ad_h\xi).$ \item $\ll C(h\cdot
v)(\Ad_{h}\eta),h\cdot w\gg_\Sl=\ll C(v)(\eta),w\gg_\Sl$.
\end{enumerate}
\end{lemma}
\begin{proof}Note that $G_{p_x}\subset G_x=H$. For (1), let $\lambda\in\g$. Using the
invariance properties of the locked inertia tensor (Lemma
\ref{lockinertia}), and the $H$-equivariance of the exponential
map at $x$, we have
\begin{align*}
\langle \Ad_{h^{-1}}^*\left[(\D\II\cdot v)(\xi)\right],\lambda
\rangle & = (\D\II\cdot v)(\xi,\Ad_{h^{-1}}\lambda)  = \frac
{d}{dt}\restr{t=0}\II(\exp_x tv)(\xi,\Ad_{h^{-1}}\lambda)\\
& \hspace{-7mm} =
 \frac {d}{dt}\restr{t=0}\II(h\exp_x tv)(\Ad_h\xi,\lambda)=  \frac {d}{dt}\restr{t=0}\II(\exp_x th\cdot
v)(\Ad_h\xi,\lambda)\\
& \hspace{-7mm} = \left\langle (\D\II\cdot (h\cdot
v))(\Ad_h\xi),\lambda\right\rangle .
\end{align*}

To prove (2), we can use the definition $$\ll
C(v)(\eta),w\gg_\Sl=\ll\nabla_{\overline{\eta}}\,\overline{v}(x),w\gg,$$ and then expand
$$\begin{array}{l}
2\ll C(h\cdot v)\left(\Ad_{h}\eta),h\cdot w\gg_\Sl  =
(\Ad_h\eta)_Q(\ll \overline{h\cdot v},\overline{h\cdot w}\gg \right)(x)\vspace{1mm}\\
+\overline{h\cdot v}\left(\ll\overline{h\cdot w},(\Ad_h\eta)_Q\gg\right)(x)
-\overline{h\cdot w}\left(\ll (\Ad_h\eta)_Q,\overline{h\cdot v}\gg\right)(x)\vspace{1mm}\\
-\ll (\Ad_h\eta)_Q(x),\left[\overline{h\cdot v},\overline{h\cdot w}\right] (x)\gg
+\ll h\cdot v,\left[\overline{h\cdot w},(\Ad_h\eta)_Q\right](x)\gg\\
+\ll h\cdot w,\left[(\Ad_h\eta)_Q,\overline{h\cdot v}\right](x)\gg .
\end{array}$$
The last three terms of this expression vanish by
Proposition~\ref{brackets}, since $\Ad_h\eta\in\lir$ due to the
invariance $G_{p_x}$-invariance of the splitting \eqref{ghr}. We
can develop the first term as
$$\begin{array}{l}
(\Ad_h\eta)_Q\left(\ll \overline{h\cdot v},\overline{h\cdot w}\gg \right)(x)=\\
=\ll\left[(\Ad_h\eta)_Q,\overline{h\cdot v}\right](x),h\cdot w\gg+\ll h\cdot v, \left[(\Ad_h\eta)_Q,\overline{h\cdot v}\right](x) \gg=0\\
=\ll\left[\eta_Q,\overline{v}\right](x),w\gg+\ll v, \left[\eta_Q,\overline{v}\right](x)\gg\\
=\eta_Q\left(\ll \overline{v},\overline{w}\gg \right)(x),
\end{array}$$
also by Proposition~\ref{brackets} and using the fact that
fundamental vector fields for an isometric action are Killing. For
the second term we have
$$\begin{array}{l}
\overline{h\cdot v}\left(\ll\overline{h\cdot w},(\Ad_h\eta)_Q\gg\right)(x)=\\
=\frac{d}{dt}\rrestr{t=0}\ll\overline{h\cdot w}(\phi([e,th\cdot v])),(\Ad_h\eta)_Q(\phi([e,th\cdot v]))\gg\vspace{1mm}\\
=\frac{d}{dt}\rrestr{t=0}\ll\overline{h\cdot w}(\phi([e,th\cdot v])),(\Ad_h\eta)_Q(h\cdot\phi([e,t v]))\gg\vspace{1mm}\\
=\frac{d}{dt}\rrestr{t=0}\ll\overline{h\cdot w}(\phi([e,th\cdot v])),h\cdot(\eta_Q(\phi([e,tv])))\gg\vspace{1mm}\\
=\frac{d}{dt}\rrestr{t=0}\ll\frac{d}{ds}\rrestr{s=0}\phi([e,th\cdot v+sh\cdot w]),h\cdot(\eta_Q(\phi([e,t v])))\gg\vspace{1mm}\\
=\frac{d}{dt}\rrestr{t=0}\ll h\cdot\frac{d}{ds}\rrestr{s=0}\phi([e,tv+s w]),h\cdot(\eta_Q(\phi([e,t v])))\gg\vspace{1mm}\\
=\frac{d}{dt}\rrestr{t=0}\ll h\cdot(\overline{w}(\phi([e,t v]))),h\cdot(\eta_Q(\phi([e,t v])))\gg\vspace{1mm}\\
=\frac{d}{dt}\rrestr{t=0}\ll
\overline{w}(\phi([e,tv])),\eta_Q(\phi([e,tv]))\gg=\overline{v}\left(\ll\overline{w},\eta_Q\gg\right)(x).
\end{array}$$
We get an analogous result for the third term. Substituting these
terms proves (2).

\end{proof}

    \section{The infinitesimal generators of a cotangent-lifted
    action}\label{secinfinitesimal}
 In this section we obtain the  $I$-representation of the infinitesimal
cotangent-lifted action. That is, for each $\xi\in\g$ we
explicitly obtain the expression
$$I\left(\xi_{T^*Q}(p_x)\right)\in T_xQ\oplus T_x^*Q\simeq (\lir\oplus\Sl
)\oplus (\lir^*\oplus\Sl^*),$$ where
$I\left(\xi_{T^*Q}(p_x)\right)=\left(T_{p_x}\tau(\xi_{T^*Q}(p_x)),K(\xi_{T^*Q}(p_x))\right)$.

Since
 the bundle projection $\tau$ is $G$-equivariant it follows that
$T_{p_x}\tau(\xi_{T^*Q}(p_x))=\xi_Q(x)$.  Using
\eqref{splittangent}, the element $\xi_Q(x)$ is represented by
$(\xi^\lir,0)\in\lir\oplus\Sl$, and so
$$I(\xi_{T^*Q}(p_x))=\left(\xi^\lir,0;K(\xi_{T^*Q}(p_x))\right).$$
For the computation of $K(\xi_{T^*Q}(p_x))$ one needs to choose a
 curve in $T^*Q$ starting at $p_x$ which locally integrates
$\xi_{T^*Q}(p_x)$. We use the curves $\hat{c}(t)=(\Exp t\xi)\cdot
p_x$ and $\hat{c}(t)=(\exp_e t\xi)\cdot p_x$, which project by
$\tau$ to the curves $c(t)=(\Exp t\xi)\cdot x$ and $c(t)=(\exp_e
t\xi)\cdot x$ respectively. By \eqref{splittangent} and \eqref{K}
we can write $K(\xi_{T^*Q}(p_x))=(\nu,\beta)\in\lir^*\oplus\Sl^*$
where $\nu$ and $\beta$ are defined by \be \langle\nu
,\lambda\rangle =
\left\langle\frac{D^\nabla_c\hat{c}(t)}{Dt}\restr{t=0},\lambda_Q(x)\right\rangle\quad\mathrm{and}\quad
\langle\beta , w\rangle  = \left\langle
\frac{D^\nabla_c\hat{c}(t)}{Dt}\restr{t=0}  ,w\right\rangle ,\ee
for every $\lambda\in\lir$ and $w\in\Sl$.

Let us recall that when  $\xi\in\h$, the curve  $\hat{c}(t)$ is  a
curve in $T^*_xQ$ passing through $p_x$ at $t=0$, and so
$K(\xi_{T^*Q}(p_x))$ is just the derivative of $\hat{c}(t)$,
regarded as an element of $T^*_xQ$. Note that also that for every
$\xi\in\g$, due to the linearity of $K$ and of the
cotangent-lifted action on fibers we have
 \be\begin{array}{lll}\label{Klineal}
K\left(\xi_{T^*Q}(p_x+p'_x)\right) & =
K\left(\xi_{T^*Q}(p_x)\right)+K\left(\xi_{T^*Q}(p'_x)\right),\\
K\left((\xi+\lambda)_{T^*Q}(p_x)\right) & =
K\left(\xi_{T^*Q}(p_x)\right)+K\left(\lambda_{T^*Q}(p_x)\right).
 \end{array}\ee

 Recall also that   under  the isomorphism~\eqref{splittangent}, a given point
$p_x$ in $T^*_xQ$ can be expressed as
$p_x=\FL(\eta_Q(x)+s)\simeq(\mu, \alpha)\in\lir^*\oplus \Sl^*$
with $\mu=\II(x)(\eta)$ and $\alpha=\ll s,\cdot\gg_\Sl$, where
$(\eta, s)\in\lir\oplus\Sl$. We will start by characterizing the
momentum and the stabilizer of $p_x$.
\begin{prop}\label{isotropyphase}
For $p_x\simeq(\mu,\alpha)\in\lir^*\oplus\Sl^*$ and $H=G_x$ the
following hold:
\begin{enumerate}
\item $\J(p_x)=\mu$, \item $G_{p_x}=H_\alpha\cap G_\mu$.
\end{enumerate}
\end{prop}
\begin{proof}
By using the formula for the momentum map for cotangent-lifted
actions, for every $\xi\in\g$,
$$\begin{array}{lll}\langle\J(p_x),\xi\rangle & = & \langle
p_x,\xi_Q(x)\rangle =\ll
\eta_Q(x)+s,\xi_Q(x)\gg=\ll\eta_Q(x),\xi_Q(x)\gg\\  & = & \II
(x)(\eta,\xi)=\langle\mu,\xi\rangle,\end{array}$$ since $s\in\Sl
=(\g\cdot x)^\perp$ which proves (1).

For (2), recall that $G_{p_x}\subset H=G_x$, and therefore
$G_{p_x}=H_{p_x}$ where by $H_{p_x}$ we mean the stabilizer of
$p_x$ under the linear action of $H$ on $T_x^*Q$ which is the
restriction to $H$ of the cotangent-lifted action of $G$. This
linear action is expressed under the isomorphism $T^*_xQ\simeq
\lir^*\oplus\Sl^*$ as
$$h\cdot (\nu,\beta)=(\Ad_{h^{-1}}^*\nu,h\cdot\beta ).$$
Therefore, $h\in H_{p_x}$ if and only if $\Ad_{h^{-1}}^*\mu=\mu
,\quad\mathrm{and}\quad h\cdot\alpha=\alpha$. That is, $h\in
G_\mu\cap H_\alpha=G_{p_x}$.
\end{proof}

\begin{thm}\label{infinitesimalaction}
Let $p_x=\FL(s+\eta_Q(x))\simeq (\mu,\alpha)$ be a point in
$T_x^*Q$. The $I$-repre\-sentation of the infinitesimal
cotangent-lifted action at $p_x$ is given by the assignment to
each element $\xi\in\g$ the element $I(\xi_{T^*Q}(p_x))$  of
$T_{p_x}(T^*Q)\simeq (\lir\oplus\Sl)\oplus (\lir^*\oplus \Sl^*)$,
expressed as
\begin{align*} I(\xi_{T^*Q}(p_x))=\left(\xi^\lir,0\, ;\, \frac
12\Proj_{\lir}\, \left[\left(\mathbf{D}\II\cdot(
\eta_Q(x)+s)\right)(\xi^\lir)-\ad^*_{\xi^\lir}\mu\right]-\ad^*_{\xi^\h}\mu
\,,\right.\\
\left.\nonumber \xi^\h\cdot\alpha-\frac 12\Proj_{\Sl}\,
\left[\left(\mathbf{D}\II\cdot
(\cdot)\right)(\xi^\lir,\eta)\right]+\ll
C(s)(\xi^\lir),\cdot\gg_\Sl\right), \end{align*} where
$\Proj_{\lir}:\g^*\rightarrow\lir^*$ and
$\Proj_\Sl:T^*_xQ\rightarrow \Sl^*$ denote the natural projections
associated to the inclusions $\lir\hookrightarrow\g$ and
$\Sl\hookrightarrow T_xQ$, and where $C(s)$ is defined in
\eqref{defC}.
\end{thm}
\begin{proof}
Let $p_x=\FL(s+\eta_Q(x))$ with $s\in\Sl$ and $\eta\in\lir$. Then
$p_x\simeq(\mu,\alpha)$ with $\J(p_x)=\mu=\II (x)(\eta)\in\lir^*$
and $\ll s,\cdot\gg_\Sl=\alpha\in\Sl^*$. This describes every
point in the fiber $T_x^*Q$ since $T^*_xQ\simeq \lir^*\oplus
\Sl^*$. Using formulae \eqref{Klineal} we have
$$\begin{array}{l}
K\left(\xi_{T^*Q}(p_x)\right)  =
K\left((\xi^\lir)_{T^*Q}(p_x)\right)+K\left((\xi^\h)_{T^*Q}(p_x)\right)=\\
\hspace{1cm}
 = K\left((\xi^\lir)_{T^*Q}(\FL
(\eta_Q(x)))\right)+K\left((\xi^\lir)_{T^*Q}(\FL
(s))\right)+K\left((\xi^\h)_{T^*Q}(p_x)\right) .
\end{array}$$
 We
now compute each term in the last equality of the above
expression. Let us start by considering $\xi=\xi^\h$. Then,
$$\hat{c}(t)=(\Exp t\xi^\h)\cdot\FL\left(s+\eta_Q(x)\right)=
\FL\left((\Exp t\xi^\h)\cdot s+(\Ad_{\Exp
(t\xi^\h)}\eta)_Q(x)\right)$$ which projects to the constant curve
$c(t)=x$. If we represent $\hat{c}(t)$ as a curve in $\lir^*\oplus
\Sl^*$ we have
$$\hat{c}(t)=(\Ad^*_{\Exp (-t\xi^\h)}\mu,(\Exp t\xi^\h)\cdot\alpha).$$

 Recalling that for curves lying in
$\mathcal{V}_{p_x}(T^*Q)$, $K$ is just the derivative along the
fiber, we get
\be\label{A}K(\xi^\h_{T^*Q}(p_x))=(-\ad^*_{\xi^\h}\mu,\xi^\h\cdot\alpha).\ee
 For the computation of
$K\left((\xi^\lir)_{T^*Q}(\FL (s))\right)$, consider
$$\hat{c}(t)=(\exp_e t\xi^\lir)\cdot \FL (s)=\FL((\exp_e
t\xi^\lir)\cdot s).$$ Then, for any $\lambda\in\lir$, \be
\label{contribution1}\left\langle\frac{D^\nabla_c\hat{c}(t)}{Dt}\restr{t=0}
,\lambda_Q(x)\right\rangle   =  \ll
\nabla_{\overline{\xi}}\,\overline{s},\lambda_Q(x)\gg =  \frac
12(\D\II\cdot s)(\xi^\lir,\lambda),\ee  by Proposition
\ref{covariant}-(3). Here we have used the fact that
$\overline{\xi}(x)=\frac{d}{dt}\rrestr{t=0}c(t)$ and that, by the
proof of Lemma \ref{orto},  $\overline{s}((\exp_et\xi^\lir)\cdot
x)=(\exp_et\xi^\lir)\cdot s$ for $t$ small, so $\FL(\overline{s})$
is a local extension of $\hat{c}(t)$. Also, for every $v\in\Sl$,
\be\label{contribution2}
\left\langle\frac{D^\nabla_c\hat{c}(t)}{Dt}\restr{t=0}
,v\right\rangle  =
\ll\nabla_{\overline{\xi}}\,\overline{s}(x),v\gg=\ll
C(s)(\xi^\lir),v\gg_\Sl \ee by \eqref{defC}. Therefore,
\begin{equation}\label{inf1}K(\xi^\lir_{T^*Q}(\FL (s)))=\left(\frac
12\Proj_\lir\,\left[(\D\II\cdot s)(\xi^\lir)\right],\ll
C(s)(\xi^\lir),\cdot\gg_\Sl\right).\end{equation}

For the computation of $K\left((\xi^\lir)_{T^*Q}(\FL
(\eta_Q(x)))\right)$, let us consider $$\hat{c}(t)=(\Exp
t\xi^\lir)\cdot \FL (\eta_Q(x))=\FL((\Exp t\xi^\lir)\cdot
\eta_Q(x)).$$ We consider also the vector field $P$ along the
local curve $c(t)=(\Exp t\xi^\lir)\cdot x$ given by $P((\Exp
t\xi^\lir)\cdot x)=(\Exp t\xi^\lir)\cdot \eta_Q(x)=(\Ad_{\Exp
(t\xi^\lir)}\eta)_Q((\Exp t\xi^\lir)\cdot x)$.
 Then, for any $v\in\Sl$,
\be\label{contribution3}\begin{array}{lll}
\left\langle\frac{D^\nabla_c\hat{c}(t)}{Dt}\restr{t=0},v\right\rangle
& = \left((\xi^\lir)_Q\ll P,\overline{v}\gg \right)(x)-\ll
\eta_Q(x),\nabla_{\overline{\xi^\lir}}\overline{v}\gg\\
& =  - \ll \eta_Q(x),\nabla_{\overline{\xi^\lir}}\overline{v}\gg =
-\frac 12 (\D\II\cdot v)(\xi^\lir,\eta),\end{array}\ee where the
last equality holds in view of Proposition \ref{covariant}-3).
Also, successive applications of Proposition \ref{covariant} and
Lemma \ref{lockinertia} give, for every $\lambda\in\lir$,
$$\begin{array}{ll}
\left\langle\frac{D^\nabla_c\hat{c}(t)}{Dt}\restr{t=0}
,\lambda_Q(x)\right\rangle &= ((\xi^\lir)_Q\ll P,\lambda_Q\gg
)(x)-\ll
\eta_Q(x),\nabla_{\overline{\xi}^\lir}\lambda_Q(x)\gg\\
& =  \frac {d}{dt}\rrestr{t=0} \ll (\Exp
t\xi^\lir)\cdot\eta_Q(x),\lambda_Q((\Exp t\xi^\lir)\cdot x)\gg\\
&\qquad -  \ll
\eta_Q(x),\nabla_{\overline{\xi^\lir}}\lambda_Q(x)\gg\\
& =  \frac {d}{dt}\rrestr{t=0}\II ((\Exp t\xi^\lir)\cdot
x)(\Ad_{\Exp (t\xi^\lir)}\eta,\lambda)\\
 & \qquad -\ll
\eta_Q(x),\nabla_{\overline{\xi^\lir}}\lambda_Q(x)\gg\\
& =  \frac {d}{dt}\rrestr{t=0}\II (x)(\eta,\Ad_{\Exp
(-t\xi^\lir)}\lambda)-\ll
\eta_Q(x),\nabla_{\overline{\xi^\lir}}\lambda_Q(x)\gg\\
& =  -\II (x)(\eta,\ad_{\xi^\lir}\lambda)-\frac 12 (\D\II\cdot
\xi^\lir)(\eta,\lambda)+\frac 12\II
(x)(\xi^\lir,\ad_\lambda\eta)\\
& = \frac 12 \left\{-\langle\ad^*_{\xi^\lir}\mu
,\lambda\rangle-\II (x)(\ad_\eta\xi^\lir,\lambda)-\II
(x)(\ad_\eta\lambda,\xi^\lir)\right\}\\
& =  \frac 12\left\{-\langle\ad^*_{\xi^\lir}\mu
,\lambda\rangle+(\D\II\cdot\eta_Q(x))(\lambda,\xi^\lir)\right\}.
\end{array}$$
Therefore, this last equation and \eqref{contribution3} yield
\begin{equation}\label{inf2}K\left(\xi^\lir_{T^*Q}(\FL(\eta_Q(x)))\right)=\left(\frac
12\Proj_{\lir}\left[(\D\II\cdot\eta_Q(x))(\xi^\lir)-\ad^*_{\xi^\lir}\mu\right],-\frac
12\Proj_\Sl\left[(\D\II\cdot (\cdot
))(\xi^\lir,\eta)\right]\right).
\end{equation}
The result of the theorem is now a consequence of \eqref{A},
\eqref{inf1} and \eqref{inf2}.
\end{proof}

\paragraph{{\bf Remark:}}In view of Theorem \ref{infinitesimalaction}, it is
immediate that $I(\xi_{T^*Q}(p_x))=0$, (and hence
$\xi_{T^*Q}(p_x)=0$) if and only if
$$\xi^\lir=0,\quad \xi^\h\cdot\alpha=0,\quad\mathrm{and}\quad
 \ad^*_{\xi^\h}\mu =0$$
hold simultaneously. This is equivalent to $\xi\in
\mathrm{Lie}(H_\alpha)\cap \mathrm{Lie}(G_\mu)$, which by
Proposition~\ref{isotropyphase} is the condition
$\xi\in\mathrm{Lie}(G_{p_x})$.

   \section{The symplectic normal space of a cotangent-lifted
action}\label{secV}
  Let $\mathcal{P}=T^*Q$ as before, be endowed
with the cotangent lift of the action of $G$ on $Q$, and let
$p_x\simeq (\mu,\alpha)$ with $\J(p_x)=\mu$. We will characterize
the symplectic normal space $N$ at $p_x$. This characterization
will be achieved by constructing an explicit choice of the
subspace $V$ appearing in the splitting \eqref{introV}. We will
continue in the setup of the previous sections.

First of all, we need a concrete choice of the $G_{p_x}$-invariant
splitting $\g=\h\oplus\lir$ introduced in \eqref{ghr}. To obtain
this we proceed as follows: let
$\h_\alpha=\mathrm{Lie}(H_\alpha),\,\g_{p_x}=\mathrm{Lie}(G_{p_x})$.
Let $H_\mu$ be the stabilizer of $\mu$ with respect to the
restriction to $H$ of the coadjoint representation. of $G$. We
note that $H_\mu=G_\mu\cap H$, and so
$\h_\mu=\mathrm{Lie}(H_\mu)=\g_\mu\cap\h$. Start by choosing a
$G_{p_x}$-invariant complement $\mathfrak{p}$ to $\h_\mu$ in
$\g_\mu$. Notice that $\h+\g_\mu=\h\oplus \mathfrak{p}$. Now let
$\OO\subset\g^*$ be the coadjoint orbit through $\mu$. Using the
infinitesimal generator map for the coadjoint representation
$\xi\mapsto \xi\cdot\mu=\ad^*_\xi\mu$, we can write
$T_\mu\OO=\g\cdot\mu$. This is a symplectic linear with the
$(-)$-Konstant-Kirillov-Souriau (KKS) symplectic form
\be\label{KKSform}
\Omega_{\mu}(\ad_{\lambda^1}\mu,\ad_{\lambda^2}\mu)=-\langle\mu,\ad_{\lambda^1}\lambda^2\rangle.\ee
The Lie group $H$ acts on $\OO$ by restriction of the transitive
$\Ad^{-1*}_G$-action. This action is Hamiltonian with momentum map
$\J_{\OO}:\OO\rightarrow \h^*$ given by
$\J_\OO(\nu)=-\Proj_\h[\nu]$. Note that, since $\mu\in\h^\circ$,
then $\J_{\OO}(\mu)=0$. Then the $H$-orbit through $\mu$ is
isotropic in $(T_\mu\OO,\Omega_\mu)$. It is a consequence of a
standard tool in symplectic linear algebra known as the Witt-Artin
decomposition (see \cite{OrRa2003}) that we can split $T_\mu\OO$
as
$$T_\mu\OO=\h\cdot\mu\oplus V_\mu\oplus W,$$ where
$V_\mu$ is a symplectic linear space complementary to $\h\cdot\mu$
in $(\h\cdot\mu)^{\Omega_\mu}$ (isomorphic to the symplectic
normal space at $\mu$ for the $H$-action) and $W$ is a Lagrangian
complement to $\h\cdot\mu$ in $V_\mu^{\Omega_\mu}$. Moreover,
these complements can be chosen to be $H_\mu$-invariant, in
particular $G_{p_x}$-invariant. Using the infinitesimal generator
map $\g\rightarrow T_\mu\OO$ let $\g_1,\g_2\subset\g$ be defined
by $V_\mu=\g_1\cdot\mu$ and $W=\g_2\cdot\mu$. Since
$G_{p_x}\subset G_\mu$, then
\be\label{propertyad}\Ad_{g^{-1}}^*(\ad^*_\lambda\mu)=\ad^*_{\Ad_g\lambda}\mu\quad
\forall\lambda\in\g,\,g\in G_{p_x}.\ee This shows that $\g_1$ and
$\g_2$ are $G_{p_x}$-invariant. Clearly $\g_\mu=\g_1\cap\g_2$.
Define $\q^\mu$ and $\mathfrak{k}$ as $G_{p_x}$-invariant
complements to $\g_\mu$ in $\g_1$ and $\g_2$ respectively. Then
the restrictions of the infinitesimal generator map
$\q^\mu\rightarrow V_\mu$ and $\mathfrak{k}\rightarrow W$ are
$G_{p_x}$-equivariant isomorphisms. It follows that
\be\label{splitg}\g=\h\oplus\mathfrak{p}\oplus\mathfrak{q}^\mu\oplus\mathfrak{k}\ee
is a $G_{p_x}$-invariant splitting and from \eqref{ghr} the
complement $\lir$ is defined as
$$\lir=\mathfrak{p}\oplus\q^\mu\oplus\mathfrak{k}.$$

Note that the infinitesimal generator map induces an equivariant
isomorphism between $\q^\mu$ and $N_\mu=\ker
T_\mu\J_\OO/\h\cdot\mu$, the symplectic normal space at $\mu$ for
the $H$-action on $\OO$. Also, by \eqref{KKSform} it follows that
\begin{eqnarray}\label{qmu}\Proj_\h[\ad_\lambda^*\mu] & = & 0,\quad\forall\lambda\in\q^\mu,\quad\mathrm{and}\\
\label{k}\Proj_{\q^\mu\oplus\mathfrak{k}}[\ad_\lambda^*\mu] & = &
0,\quad\forall\lambda\in\mathfrak{k}, \end{eqnarray} where the
second property uses the facts that $W$ is Lagrangian and that
$W\subset V_\mu^{\Omega_\mu}$.

 We now introduce some notation. If
$A,B$ are linear subspaces of $\Sl$ and $\Sl^*$ respectively, we
denote by $A^\circ_{\Sl^*}$ and $B^\circ_\Sl$ its annihilators in
$\Sl^*$ and $\Sl$. Similarly, with respect to the induced inner
product $\ll\cdot,\cdot\gg_\Sl$ in $\Sl$ (and the corresponding
one in $\Sl^*$, $\ll\cdot,\cdot\gg_{\Sl^*}$), $A^\perp_\Sl$ and
$B^\perp_{\Sl^*}$ denote the orthogonal of $A$ in $\Sl$ and of $B$
in $\Sl^*$. The dual space $A^*$ is identified with a subspace of
$\Sl^*$ by $A^*=\ll A,\cdot\gg_\Sl$. In particular, for any
subspace $\g'\subset\h$ we have
\be\label{dualsann}[\g'\cdot\alpha]_\Sl^\circ=(\g'\cdot
s)_\Sl^\perp\quad\text{and}\quad
([\g'\cdot\alpha]_\Sl^\circ)^*=(\g'\cdot\alpha)_{\Sl^*}^\perp .\ee
\begin{lemma}\label{lemmaaux} Let $(\h_\mu\cdot
s)^\perp_\h$ be the orthogonal complement to $\h_\mu\cdot s$ in
$\h\cdot s$ with respect to $\ll\cdot,\cdot\gg_\Sl$. Then
\be\label{isobasecot}[\h_\mu\cdot\alpha]_\Sl^\circ=(\h_\mu\cdot
s)^\perp_\h\oplus[\h\cdot \alpha]_\Sl^\circ,\ee
\end{lemma}
\begin{proof}
Use $\ll\cdot,\cdot\gg_\Sl$ to obtain the orthogonal splittings
\begin{eqnarray*}\h\cdot s=\h_\mu\cdot s\oplus (\h_\mu\cdot
s)^\perp_\h\quad\mathrm{and}\quad\Sl=\h\cdot s\oplus(\h\cdot
s)^\perp=\h_\mu\cdot s\oplus (\h_\mu\cdot
s)^\perp_\h\oplus(\h\cdot s)_\Sl^\perp .\end{eqnarray*} From
\eqref{dualsann} the result follows.
\end{proof}
Notice that by construction, all the spaces involved in
\eqref{isobasecot} are $G_{p_x}$-invariant. We define
$\mathrm{pr}_1:[\h_\mu\cdot\alpha]_\Sl^\circ\rightarrow(\h_\mu\cdot
s)^\perp_\h$ to be the equivariant projection onto the first
component of the splitting \eqref{isobasecot}.

We introduce for later use the diamond notation for a linear
representation.
\begin{defn}Given a linear space $L$ supporting a
representation of a compact Lie group $M$ with Lie algebra
$\mathfrak{m}$ and $(l,o)\in L\times L^*$, let
$l\diamond_\frak{m}o\in\frak{m}^*$ be defined as
$$\langle l\diamond_\frak{m} o,\xi\rangle=\langle o,\xi\cdot
l\rangle,$$ for any $\xi\in \frak{m}$.
\end{defn}

The next theorem provides the first main result of the paper: an
explicit choice of the symplectic normal space for a
cotangent-lifted action as a subspace of $T_{p_x}(T^*Q)$, this
linear space being identified with
$(\lir\oplus\Sl)\oplus(\lir^*\oplus\Sl^*)$ through the
$I$-representation.

\begin{thm}\label{generalsymplectic}
The symplectic normal space $V$ at the point $p_x\simeq
(\mu,\alpha)=\FL(\eta,s)$, with $\eta\in\lir,\,s\in\Sl$, and
$G_x=H$ embeds linearly into $T_{p_x}(T^*Q)$ as the following
$G_{p_x}$-invariant subspace: \be\label{Vgeneral} V  =
\mathrm{span}\left\langle\left(
\lambda+j\left(\mathrm{pr}_1(a)\right),\,a\,  ;\,
f_1(\lambda+j(\mathrm{pr}_1(a)),a)
,\,\beta+f_2(\lambda+j(\mathrm{pr}_1(a)),a)\,\right)\right\rangle,
\ee where $$  \begin{array}{rl}
\lambda\in & \q^\mu,\\
a\in  & [\h_\mu\cdot\alpha]_\Sl^\circ,\\
\beta\in &
([\h_\mu\cdot \alpha]_\Sl^\circ)^*,\vspace{2mm}\\
f_1(\gamma)= & \frac 12 \Proj_\lir\, \left[ \left(\D\II \cdot
(\eta_Q(x)+s)\right)(\gamma) -(\D\II \cdot a)(\eta)  +
\ad^*_\gamma\mu
\right]\\ & +\ll C(s)(\cdot),a\gg_\Sl ,\vspace{2mm}\\
f_2(\gamma)= & -\frac 12\Proj_\Sl\left[\left(\D\II\cdot
(\cdot)\right)(\eta,\gamma)\right]+\ll C(s) (\gamma),\cdot\gg_\Sl.
\end{array}$$

The spaces $\q^\mu$ and $\mathfrak{k}$ are components of the
splitting \eqref{splitg}. The linear map $j:(\h_\mu\cdot
s)^\perp_\h\rightarrow \mathfrak{k}$ is defined by $\Proj_\h\,
[\ad^*_{j(b)}\mu]-b\diamond_\h\alpha=0$, for all $b\in
(\h_\mu\cdot s)^\perp_\h$.
\end{thm}
\begin{rem}
Of course it is clear that many other realizations of the space
$V$ can be constructed. However the one introduced in
Theorem~\ref{generalsymplectic} has the advantage of putting its
inherited symplectic form into normal form (see
Corollary~\ref{generalsymplectic2}).
\end{rem}
\begin{proof}
It follows from \eqref{sasakisymplectic} and
Theorem~\ref{infinitesimalaction} that an element
$(\lambda,a;\,\nu,\beta)\in T_{p_x}(T^*Q)$ is symplectically
orthogonal to $\g\cdot p_x$ if and only if
\begin{align*} \langle\nu ,\xi^\lir\rangle - \frac
12\left\{(\D\II\cdot(\eta_Q(x)+s))(\xi^\lir,\lambda)-(\D\II\cdot
a)(\xi^\lir,\eta)-\langle\ad^*_{\xi^\lir}\mu ,\lambda\rangle\right\}\\
- \ll C(s)(\xi^\lir),a\gg_\Sl+\langle\ad^*_{\xi^\h}\mu
,\lambda\rangle-\langle \xi^\h\cdot\alpha,a\rangle=0,
\end{align*}
 for every $\xi^\lir\in\lir$ and $\xi^\h\in\h$.
Let us define the map $R:\lir\oplus\Sl\rightarrow\h^*$ by
$$ R(\lambda,a)  =  \Proj_\h\,
[\ad^*_\lambda\mu]-a\diamond_\h\alpha.$$ Define the vector
subspaces
\begin{eqnarray*}V' & = &
\left\{(\lambda,a;\nu,\beta)\in(\lir\oplus\Sl)\oplus(\lir^*\oplus\Sl^*)\,:\,\nu=f_1(\lambda,a)\right\},\quad
\mathrm{and}\\
\overline{R} & = &
\left\{(\lambda,a;\nu,\beta)\in(\lir\oplus\Sl)\oplus(\lir^*\oplus\Sl^*)\,:\,(\lambda,a)\in
\ker R\right\}.\end{eqnarray*} It is clear that $\left(\g\cdot
p_x\right)^\w =V'\cap\overline{R}$.  Recall from \eqref{introV}
that $V$ is a ($G_{p_x}$-invariant) complement to $\g_\mu\cdot
p_x$ in $(\g\cdot p_x)^\w$. From Theorem \ref{infinitesimalaction}
we obtain \be\label{gmuexpression}\begin{array}{lll}\g_\mu\cdot
p_x & = \left\{\left(\xi^\lir,0;\,\frac 12 \Proj_{\lir}\,
\left[\left(\mathbf{D}\II\cdot(
\eta_Q(x)+s)\right)(\xi^\lir)\right],\right.\right.\\
&\left.\left. \xi^\h\cdot\alpha-\frac
12\Proj_{\Sl}\,\left[\left(\mathbf{D}\II\cdot (\cdot
)\right)(\xi^\lir,\eta)\right]+\ll
C(s)(\xi^\lir),\cdot\gg_\Sl\right)\,
:\,\forall\,\xi^\lir\in\mathfrak{p},\,\xi^\h\in\h_\mu\right\}.\\
\end{array}\ee
Since $\g_\mu\cdot p_x\subset V'\cap\overline{R}$, if we find a
$G_{p_x}$-invariant splitting $V'=\g_\mu\cdot p_x\oplus V''$ then
$V=V''\cap\overline{R}$.  From the above expressions, it is clear
that we can choose $V''$ as
$$V''=\left\{(\lambda,a;f_1(\lambda,a),\beta+f_2(\lambda)\,)\,:
\,\lambda\in\q^\mu\oplus\mathfrak{k},\,a\in\Sl,\,\beta\in
([\h_\mu\cdot \alpha]_\Sl^\circ)^*\right\}.$$ Finally, to obtain
$V$, recall from the proof of Lemma \ref{lemmaaux} that
$$\Sl=\h_\mu\cdot s\oplus (\h_\mu\cdot s)^\perp_\h\oplus [\h\cdot
\alpha]_\Sl^\circ .$$ Then we can then write
$\lambda=\lambda_1+\lambda_2$ and $a=a_1+a_2+a_3$ with
$\lambda_1\in\q^\mu,\,\lambda_2\in\mathfrak{k},\,a_1\in
\h_\mu\cdot s,\,a_2\in (\h_\mu\cdot s)^\perp_\h$ and $a_3\in
[\h\cdot\alpha]_\Sl^\circ$. By studying $R$ and its composition
with the projection $\h^*\rightarrow\h_\mu^*$, it follows that
$R(\lambda, a)=0$ if and only if $a_1=0$ and $R(\lambda_2,a_2)=0$.
Let $R':\mathfrak{k}\oplus (\h_\mu\cdot s)^\perp_\h\rightarrow
\h^*$ be the restriction $R\rrestr{\mathfrak{k}\oplus (\h_\mu\cdot
s)^\perp_\h}$. Therefore, $V$ is characterized by
\be\label{intermediateV} V
=\mathrm{span}\,\langle(\lambda+\gamma,a+b;f_1(\lambda+\gamma,a+b),\beta+f_2(\lambda+\gamma)\,)\rangle,\ee
with $\lambda\in\q^\mu,\,a\in [\h\cdot \alpha]_\Sl^\circ,\,
(\gamma,b)\in \ker R'$ and $\beta\in ([\h_\mu\cdot
\alpha]_\Sl^\circ)^*$.

Next, notice that there is no $0\neq\gamma\in\mathfrak{k}$ such
that $(\gamma,0)\in \ker R'$, since this amounts to $\Proj_\h
[\ad^*_\gamma\mu]=0$, which is a contradiction with
$\gamma\in\mathfrak{k}$. Hence, there is a linear subspace
$D\subset (\h_\mu\cdot s)^\perp_\h$ and a linear
$G_{p_x}$-equivariant map $j:D\rightarrow\mathfrak{k}$ such that
\be\label{kerR'}\ker R'=\left\{(j(b),b)\, :\, b\in D\right\}.\ee
The map $j$ induces an isomorphism between $\ker R'$ and $D$, and
in particular it follows that $\dim \ker R'=\dim D$.

We now prove $D=(\h_\mu\cdot s)^\perp_\h$.  From
\eqref{intermediateV}, it follows that
$$\dim D = \dim V-\dim \q^\mu-\dim [\h\cdot\alpha]^\circ_\Sl-\dim
[\h_\mu\cdot\alpha]^\circ_\Sl.$$ Note that $V$ and $\q^\mu$ are
symplectic normal spaces identified with $\ker
T_{p_x}\J/\g_\mu\cdot p_x$ and $\ker T_\mu\J_{\OO}/\h\cdot\mu$. By
the Bifurcation Lemma (see \cite{OrRa2003}), given a Hamiltonian
action of a Lie group $G$ on $\PP$, and $z\in\PP$, the relation
$\im T_{z}\J=[\mathfrak{g}_z]^\circ\subset\mathfrak{g}^*$ holds.
This implies that \begin{eqnarray*} \dim V & = & \dim T^*Q-(\dim
G-\dim G_{p_x})-\dim G_\mu\cdot p_x\\ & = & 2\dim\Sl+\dim G+2\dim
G_{p_x}-2\dim H-\dim G_\mu\\
\dim \q^\mu & = & \dim \OO -(\dim H-\dim H_\mu)-\dim H\cdot\mu\\
& = & \dim G+2\dim H_\mu-\dim G_\mu-2\dim H
\end{eqnarray*}
Note that \eqref{isobasecot} implies $\dim
[\h\cdot\alpha]^\circ_\Sl+\dim [\h_\mu\cdot\alpha]^\circ_\Sl=2\dim
[\h_\mu\cdot\alpha]^\circ_\Sl-\dim (\h_\mu\cdot s)_\h^\perp$.
Putting all contributions together we have
\begin{eqnarray*}
\dim D & = & 2\dim \Sl+2\dim G_{p_x}-2\dim H_\mu-2\dim
[\h_\mu\cdot\alpha]^\circ_\Sl + \dim
(\h_\mu\cdot s)_\h^\perp\\
& = & 2\dim \Sl+2\dim G_{p_x}-2\dim H_\mu -2(\dim\Sl-\dim
H_\mu+\dim G_{p_x})+ \dim (\h_\mu\cdot s)_\h^\perp\\ & = & \dim
(\h_\mu\cdot s)_\h^\perp,\end{eqnarray*} where we have used that
${(H_\mu)}_\alpha=H_\mu\cap H_\alpha=G_{p_x}$. Therefore, by
linearity, $D=(\h_\mu\cdot s)^\perp_\h$, and the map $j$ is a
$G_{p_x}$-equivariant linear embedding $j:(\h_\mu\cdot
s)^\perp_\h\rightarrow \mathfrak{k}$. By Lemma \ref{lemmaaux} for
any $a\in [\h\cdot \alpha]_\Sl^\circ$ and $b\in(\h_\mu\cdot
s)^\perp_\h$ there is a unique $a'\in
[\h_\mu\cdot\alpha]_\Sl^\circ$ such that $\mathrm{pr_1}(a')=b$ and
$a'=b+a$. This fact, together with \eqref{intermediateV} implies
that
\begin{eqnarray*} V
=\mathrm{span}\,\langle\left(\lambda+j(\mathrm{pr_1}(a')),a';f_1(\lambda+j(\mathrm{pr_1}(a')),a'),\beta+f_2(\lambda+j(\mathrm{pr_1}(a')))\,\right)\rangle,\end{eqnarray*}
with $\lambda\in\q^\mu,\,a'\in [\h_\mu\cdot \alpha]_\Sl^\circ$ and
$\beta\in ([\h_\mu\cdot \alpha]_\Sl^\circ)^*$, as stated in the
theorem.

Let us now show that $V$ is $G_{p_x}$-invariant. By construction,
the spaces $\q^\mu, [\h_\mu\cdot\alpha]^\circ$ and
$([\h_\mu\cdot\alpha]^\circ)^*$ are invariant. Then it suffices to
prove the $G_{p_x}$-equivariance of $f$.

Let $g\in G_{p_x}$. By $(1)$ in Lemma \ref{intermch3}, we easily
obtain
\begin{equation}\label{no1}(\D\II\cdot(g\cdot a))(\eta)=(\D\II\cdot(g\cdot
a))(\Ad_g\eta)=\Ad_{g^{-1}}^*((\D\II\cdot a)(\eta))\end{equation}
and \begin{equation}\label{no2}(\D\II\cdot
(\eta_Q(x)+s))(\Ad_g\lambda)=(\D\II\cdot
(g\cdot(\eta_Q(x)+s)))(\Ad_g\lambda)=\Ad^*_{g^{-1}}\left[(\D\II\cdot(\eta_Q(x)+s))(\lambda)\right].\end{equation}
By $(2)$ in Lemma \ref{intermch3} for any $\xi\in\lir$, $$\ll
C(s)(\xi),g\cdot a\gg=\ll C(g^{-1}\cdot s)(\Ad_{g^{-1}}\xi),a\gg
=\ll C(s)(\Ad_{g^{-1}}\xi),a\gg,$$ and so
\begin{equation}\label{no3}\ll C(s)(\cdot),g\cdot
a\gg=\Ad^*_{g^{-1}}(\ll C(s)(\cdot), a\gg).\end{equation} Finally,
the equivariance of $f$ follows from \eqref{propertyad},
\eqref{no1}, \eqref{no2}, \eqref{no3} and from the equivariance of
the projection $\Proj_\lir$, which is in turn a consequence of the
invariance of the splitting $\g=\h\oplus\lir$.
\end{proof}

 With the help of Theorem~\ref{generalsymplectic}
we can provide a characterization of the symplectic normal space
for a cotangent-lifted action, together with its symplectic form
and momentum map which depends solely on the coadjoint
representation of $G$ and on its isometric action on the base $Q$.

\begin{cor}\label{generalsymplectic2} Let $p_x\simeq (\mu,\alpha)=\FL
(\eta,s)$, $H=G_x$ and $B=[\h_\mu\cdot\alpha]^\circ_\Sl$. Let
$N_\mu$ be the symplectic normal space at $\mu$ for the restricted
action of $H$ on $\OO$.\\ Then
 the symplectic normal space $N$ at $p_x$
 is $G_{p_x}$-equivariantly
symplectomorphic to $N_\mu\oplus T^*B$ with symplectic form
\be\begin{array}{ccccc}
 &  & N_\mu & B & B^*\\ &&&&\\
\Omega & = & \left(
 \begin{array}{c} \Xi \\ 0 \\ 0 \end{array}\right. & \left.
\begin{array}{c} 0 \\
0\\ -\mathbf{1} \end{array} \right. & \left. \begin{array}{c} 0 \\
\mathbf{1}
\\ 0 \end{array}\right).
\end{array}\ee
That is, $\Omega=\Xi+\Omega_B$, where the symplectic form $\Xi$ on
$N_\mu$ is defined by
$$\Xi(\ad^*_{\lambda^1}\mu,\ad^*_{\lambda^2}\mu )
 =
-\langle\mu,\ad_{\lambda^1}\lambda^2\rangle$$ and $\Omega_B$ is
the canonical symplectic form on $T^*B$.
 The action of $G_{p_x}=H_\alpha\cap G_\mu$ on
$N$ is given by the expression
$$h\cdot \left(\ad^*_\lambda\mu;(a,\beta)\right)
=\left(\ad^*_{\Ad_h\lambda}\mu;(h\cdot a,h\cdot \beta)\right) ,$$
for $\ad^*_\lambda\mu\in N_\mu$, and $(a,\beta)\in T^*B$, where
$h\cdot\beta$ refers to the contragredient representation of $H$
on $\Sl^*$. The corresponding momentum map $\J_N$ is \be
\J_N\left(\ad^*_\lambda\mu;(a,\beta)\right)=\frac 12 \lambda
\diamond_{\g_{p_x}}\ad^*_\lambda\mu +a\diamond_{\g_{p_x}}\beta
.\ee
  Moreover, $N$
embeds linearly and $G_{p_x}$-equivariantly into $T_{p_x}(T^*Q)$
by the map $\iota_N:N\rightarrow
(\lir\oplus\Sl)\oplus(\lir^*\oplus\Sl^*)\simeq T_{p_x}(T^*Q)$
given by
\be\label{embeddingN}\iota_N(\ad^*_\lambda\mu;(a,\beta))=\left(\tilde{\lambda}+j(\mathrm{pr_1}(a)),a;
f_1(\tilde{\lambda}+j(\mathrm{pr_1}(a)),a),
\beta+f_2(\tilde{\lambda}+j(\mathrm{pr_1}(a)))\,\right),\ee where
$\tilde{\lambda}\in \q^\mu$ is the unique element in $\q^\mu$ such
that $\ad_{\tilde{\lambda}}^*\mu=\ad^*_\lambda\mu$, and $f_1,f_2$
are defined in the statement of Theorem~\ref{generalsymplectic}.
The embedding $\iota_N$ is a symplectomorphism onto its image $V$,
equipped with the restriction of the symplectic form of
$T_{p_x}(T^*Q)$.
\end{cor}
\begin{proof}
 It easily follows from
Theorem~\ref{generalsymplectic} that $\iota_N$ maps $N$
isomorphically and $G_{p_x}$-equivariantly to $V$. The expression
for $\Omega$ follows from its definition
$\Omega=i_N^*(\w\rrestr{V})$, and using \eqref{sasakisymplectic},
and \eqref{k}. Therefore $\iota_N$ is symplectic. The momentum map
for a symplectic linear action of a group on a symplectic linear
space $(V,\Omega)$ is defined by $\langle\J_V(v),\xi\rangle= \frac
12\Omega(\xi\cdot v,v)$. This together with \eqref{k} and noting
that $\mathfrak{k}$ is $G_{p_x}$-invariant gives the expression
for $\J_N$.
\end{proof}

\paragraph{{\bf Remarks:}} \begin{enumerate}\item While the embedding $\iota_N$ depends on the choice of an invariant metric on $Q$ and
the splitting of $\g$, the characterization of
$N$ is completely general. \item Note that the restriction $\Xi$
of $\Omega$ to $N_\mu$ is precisely the symplectic form inherited
from the KKS
form on $\OO$.\\
\item The embedding $\iota_N$ of
Corollary~\ref{generalsymplectic2}, depending on the map $j$ is
not explicit, since $j:\dim
(\h_\mu\cdot\alpha)_\h^\perp\rightarrow \mathfrak{k}$ is defined
through the kernel of a linear map. However, in a variety of
relevant situations the map $j$ is trivial, and $\iota_N$ is
totally explicit. These cases are studied in the next section.
\end{enumerate}

    \section{Particular cases of the symplectic normal space}
\label{secspecialV}
In this section we will focus on some
particular cases for which the symplectic normal space is
simplified by the fact that the map $j$ is trivial. These cases
are justified by their geometric or dynamical interest, and we
will briefly explore the possibilities that our characterization
of the symplectic normal space of a cotangent-lifted action can
offer in current and future research.

\begin{prop}
Let $p_x\simeq (\mu,\alpha)$. If either
\begin{itemize}
\item[a)] $G_\mu=G$, or \item[b)] $\alpha=0$, or \item[c)] $G$
acts locally freely at $x$, or \item[d)] $Q$ is a manifold of
constant orbit type $(H)$, i.e. every point in $Q$ has stabilizer
conjugated to $H$, or \item[e)] $H\subset G_\mu$,
\end{itemize}
then $j=0$ and so
$[\h_\mu\cdot\alpha]_\Sl^\circ=[\h\cdot\alpha]_\Sl^\circ$.
\end{prop}
\begin{proof}
By Lemma~\ref{lemmaaux}, and since $j$ is injective,  $j=0$ if and
only if $(\h_\mu\cdot s)_\h^\perp=0$ which implies
$[\h_\mu\cdot\alpha]^\circ=[\h\cdot\alpha]^\circ$.\\  a) If
$\g=\g_\mu$, then $\h_\mu=\h$. As a consequence $\h_\mu\cdot
s=\h\cdot s$, and then $(\h_\mu\cdot s)_\h^\perp=0$.\\ b) If
$\alpha=0$ then $s=0$ and then $\h\cdot s=0$ which implies $(\h_\mu\cdot s)_\h^\perp=0$. \\
c) If $G$ acts locally
trivially at $x$ then $\h=0$ and the result follows as in $b)$.\\
d) In this case we also have $\h\cdot s=0$ since if $Q$ has
constant orbit type, the linear slice $\Sl$ is a fixed-point space
for the linear $H$-action on it. The result then follows as in c). \\
e) If $H\subset G_\mu$ then $H_\mu=H$ and the result follows as in
a).
\end{proof}

    \subsection{Case a). Totally isotropic momentum
($G_\mu=G$)}\label{totallyisotropic} The first case to consider
will be that of points of the form $p_x\simeq(\mu,\alpha)$
satisfying $G_\mu=G$, that is, the elements such that their
momentum value, $\J(p_x)=\mu$, is totally isotropic, and then
$G_{p_x}=H_\alpha$. This happens for instance at any point in
$T^*Q$ if $G$ is Abelian. Near the orbit of a point $p_x$ with
totally isotropic momentum value $\mu$, the reduced space
$\PP_\mu=\J^{-1}(\mu)/G$ is modelled on the orbit space
$\J_N^{-1}(0)/G_{p_x}$, where $N$ is the symplectic normal space
at $p_x$ on which the compact group $G_{p_x}$ acts linearly with
momentum map $\J_N$. The resultant splitting of the symplectic
normal space should give a geometrical insight into the local
properties of the topological bundle structure of $\PP_\mu$ over
$Q/G$.

 In this case, since $\g_\mu=\g$ we obtain
$\q^\mu=0$ and $B=[\h\cdot\alpha]^\circ$. From Corollary
~\ref{generalsymplectic2} we get that $N$ is $H_\alpha$-isomorphic
to $T^*B$ equipped with the symplectic form \be\begin{array}{cccc}
 & & B & B^*\\
\Omega & = & \left(\begin{array}{c} 0 \\
-\mathbf{1}\end{array}\right. & \left.\begin{array}{c} \mathbf{1} \\
0\end{array} \right).
\end{array}\ee
 That is, the symplectic normal space is a cotangent bundle with
its canonical symplectic form. The $H_\alpha$-action on $T^*B$ is
by diagonal (i.e. cotangent-lifted) action and its associated
momentum map $\J_N:T^*B\rightarrow \h_\alpha^*$ has the expression
\be \J_N (a,\beta)=a\diamond_{\h_\alpha}\beta .\ee

    \subsection{Case b). Vertical covectors
($\alpha=0$)}\label{relativeequilibria} The second important case
is that of covectors which are vertical for the group action, i.e.
when $\alpha=0$ and consequently $p_x\simeq(\mu,0)$. These are the
points which are candidates to be relative equilibria in a
symmetric simple mechanical system with kinetic energy given by
the Riemannian structure in $Q$. Indeed, every relative
equilibrium of a simple mechanical system with momentum $\mu$ is
of this form, i.e. $p_x=\FL(\eta_Q(x))$, with $\eta\in\g$
satisfying the relation $\II (x)(\eta)=\mu$ (see \cite{AbMa,
MarLec} for details). We call these points vertical covectors. The
splitting of the symplectic normal space at this class of points
obtained below has important consequences in the study of the
dynamics of relative equilibria for simple mechanical systems.
Those are aspects beyond the scope of this work, but in
\cite{Miguelstability} some of the results in this section are
applied to test the orbital stability of relative equilibria at
singular values of the momentum map, generalizing the
constructions of \cite{SiLeMar} in the regular case. From
Corollary \ref{generalsymplectic2} we obtain that in this case the
symplectic normal space $N$ is $G_{p_x}$-equivariantly
symplectomorphic to $N_\mu\oplus T^*\Sl$, with symplectic form
\be\label{symplecticmatrixRE}\begin{array}{ccccc}
 & & N_\mu & \Sl & \Sl^*\\
\Omega & = & \left(\begin{array}{c} \Xi \\ 0 \\ 0
\end{array}\right. &
\begin{array}{c} 0 \\
0\\ -\mathbf{1} \end{array} & \left. \begin{array}{c} 0 \\
\mathbf{1}
\\ 0 \end{array}\right)
\end{array}\ee
where, as before,$$ \Xi(\ad^*_{\lambda_1}\mu,\ad^*_{\lambda_2}\mu)
= -\langle\mu,\ad_{\lambda_1}\lambda_2\rangle .$$

The $G_{p_x}$-action on $N_\mu \oplus T^*\Sl$ is diagonal with
associated momentum map \be\label{momentummapRE}
\J_N(\ad^*_\lambda\mu;(a,\beta))=\frac 12 \lambda\diamond_{\h_\mu}
\ad^*_{\lambda}\mu +a\diamond_{\h_\mu} \beta.\ee

      \subsection{Case c). $\h=0$.}
Much work has been done for this case from both the global and
local points of view of reduction and its applications to
geometric mechanics. In the particular case when the action of $G$
on $Q$ is globally free the quotient $Q/G$ is a manifold. It is
well known, (see \cite{MarLec,PerMar}) that under this assumption
the reduced space $\J^{-1}(\mu)/G_\mu$ can be realized as a bundle
over $T^*(Q/G)$ having as typical fiber the coadjoint orbit $\OO$.
Since in a free action situation the symplectic normal space at
$p_x$ is isomorphic to the tangent space to $\J^{-1}(\mu)/G_\mu$
at $[p_x]$, we expect to obtain that $N\simeq T_\mu\OO\oplus
T_{[p_x]}(T^*(Q/G))$. We show now that when $\h=0$ this is exactly
the content of Corollary \ref{generalsymplectic2}.

Let $p_x\simeq(\mu,\alpha)$ be a point in $T^*Q$ such that
$\g_x=\h=0$. Then $B=\Sl$ and $N_\mu\simeq T_\mu\OO$. It follows
that the symplectic normal space $N$ at $p_x$ is symplectomorphic
to $T_\mu\OO\oplus T^*\Sl$ with the symplectic form given by
$$\Omega_\mu+\Omega_\Sl,$$
where $\Omega_\mu$ is the KKS structure defined in \eqref{KKSform}
and $\Omega_\Sl$ the canonical symplectic form on $T^*\Sl$. If the
action of $G$ is free everywhere then $N$ is isomorphic to
$T_\mu\OO\oplus T_{[p_x]}(T^*(Q/G))$ as follows by the chain of
isomorphisms
$$T_{[p_x]}(T^*(Q/G))\simeq T^*(T_{[x]}(Q/G))\simeq T^*\Sl.$$
Since $\h=0$ (hence $\g_{p_x}=0$) the momentum map is trivial.

\subsection{Case d). $Q$ is of constant orbit type.}
Let $H=G_x$ and suppose that, for every $x'\in Q$, the group
$G_{x'}$ is conjugate to $H$. Then the quotient $Q/G$ is still a
smooth manifold, but in general the orbit map $Q\rightarrow Q/G$
does not define a principal bundle and the standard results for
regular cotangent bundle reduction do not apply. The study of
cotangent bundle reduction over a manifold of constant orbit type
is the natural step towards singular cotangent bundle reduction
after truly regular (free) reduction. The results obtained in this
situation, usually called Single Orbit Type theorems, have proved
to be useful in the fully singular generalization of cotangent
bundle reduction. For instance, in \cite{EmRo} and \cite{Schmah}
Single Orbit Type theorems for zero and totally isotropic momentum
were proved. This made possible in \cite{PerRoSD} to generalize
regular cotangent bundle reduction to the singular case when
$\mu=0$ in presence of several orbit types in $Q$. Also, for
general $\mu$, a Single Orbit Type theorem has been obtained in
\cite{Hoch} by using the so-called Weinstein representation of
gauged reduction (see \cite{PerRaprep}). A similar result based on
the alternative Sternberg representation can be found in
\cite{PerRo}. It is expected that these results, without being a
final answer in their own right, will be useful to establish a
fully singular picture of cotangent bundle reduction at arbitrary
momentum values.

Recall that in the constant orbit type case $\Sl$ and $\Sl^*$ are
fixed-point sets for the linear $H$-action, and therefore
$[\h_\mu\cdot\alpha]^\circ=\Sl$. Consequently, it follows from
Corollary \ref{generalsymplectic2} that the symplectic normal
space $N$ at $p_x\simeq (\mu,\alpha)$ is
$G_{p_x}$-symplectomorphic to $N_\mu\oplus T^*\Sl$ equipped with
the symplectic form given in \eqref{symplecticmatrixRE} and
momentum map given by \be \J_N(\ad^*_\lambda\mu;(a,\beta))=\frac
12\lambda\diamond_{\h_\mu} \ad^*_{\lambda}\mu.\ee
\subsection{Case e). $H\subset G_\mu$.}
In Lemma 4.1 of \cite{Schmah}, it is proved that if $H$ is a
normal subgroup of $G$, then $H\subset G_\mu$. This justifies the
study of the rather more general situation $H\subset G_\mu$ as a
particular case of Corollary \ref{generalsymplectic2}. In this
case we see that $\h_\mu\cdot\alpha=\h\cdot\alpha$ and
$G_{p_x}=H_\alpha$. Besides, since $\h\cdot\mu=0$ then
$N_\mu\simeq T_\mu\OO$, as for case c). Therefore
$B=[\h\cdot\alpha]^\circ$  and $N\simeq T_\mu\OO\oplus T^*B$,
where $N$ is equipped with the symplectic form
$\Omega=\omega_\mu+\Omega_B$. Again, $\Omega_B$ denotes the
canonical symplectic form on $T^*B$.
 The momentum map $\J_N$ for the
$H_\alpha$-action is readily verified to be
 \be \J_N(\ad_\lambda^*\mu;(a,\beta))=\frac 12
\lambda\diamond_{\h_\alpha} \ad^*_{\lambda}\mu
+a\diamond_{\h_\alpha} \beta.\ee

    \end{document}